\documentclass{amsart}
\usepackage{amsthm,amsmath,amssymb,amscd,graphics,pb-diagram,enumerate}

\newcommand{\thisdate}{\today}
{
   \newtheorem{theorem}{Theorem}[subsection]                     
   \newtheorem{proposition}[theorem]{Proposition}     
   \newtheorem{lemma}[theorem]{Lemma}
   \newtheorem*{claim}{Claim}
   \newtheorem{corollary}[theorem]{Corollary}

}
{\theoremstyle{definition}

   \newtheorem{definition}[theorem]{Definition}
   \newtheorem{remark}[theorem]{Remark}
}

\newcommand{\bbS}{{\mathbb{S}}}

\newcommand{\CC}{{\mathbb{C}}}

\newcommand{\PP}{{\mathbb{P}}}
\newcommand{\ZZ}{{\mathbb{Z}}}

\newcommand{\bK}{{\mathbf{K}}}
\newcommand{\bM}{{\mathbf{M}}}
\newcommand{\bX}{{\mathbf{X}}}

\newcommand{\bmu}{{\boldsymbol{\mu}}}

\newcommand{\g}{{\mathfrak{g}}}
\newcommand{\m}{{\mathfrak{m}}}

\newcommand{\cA}{{\mathcal A}}
\newcommand{\cB}{{\mathcal B}}
\newcommand{\cC}{{\mathcal C}}
\newcommand{\cD}{{\mathcal D}}
\newcommand{\cE}{{\mathcal E}}
\newcommand{\cF}{{\mathcal F}}
\newcommand{\cG}{{\mathcal G}}
\newcommand{\cH}{{\mathcal H}}
\newcommand{\cI}{{\mathcal I}}
\newcommand{\cK}{{\mathcal K}}
\newcommand{\cL}{{\mathcal L}}
\newcommand{\cM}{{\mathcal M}}
\newcommand{\cN}{{\mathcal N}}
\newcommand{\cO}{{\mathcal O}}

\newcommand{\cR}{{\mathcal R}}
\newcommand{\cS}{{\mathcal S}}
\newcommand{\cT}{{\mathcal T}}

\newcommand{\cX}{{\mathcal X}}
\newcommand{\cY}{{\mathcal Y}}

\newcommand{\Av}{\operatorname{Av}}
\newcommand{\Id}{\operatorname{Id}}
\newcommand{\Spec}{\operatorname{Spec}}
\newcommand{\graph}{\operatorname{Graph}}
\newcommand{\Proj}{\operatorname{Proj}}
\newcommand{\Sect}{\operatorname{Sect}}
\newcommand{\Isom}{\operatorname{Isom}}
\newcommand{\Sing}{\operatorname{Sing}}

\newcommand{\Stab}{\operatorname{Stab}}

\newcommand{\Hom}{{\operatorname{Hom}}}

\newcommand{\Ext}{{\operatorname{Ext}}}
\newcommand{\cursext}{{\mathcal E}xt}

\newcommand{\Aut}{{\operatorname{Aut}}}

\newcommand{\sh}{{\operatorname{sh}}}

\newcommand{\opH}{{\operatorname{H}}}

\newcommand{\depth}{\operatorname{depth}}
\newcommand{\chara}{\operatorname{char}}

\newcommand{\lrar}{\longrightarrow}
\newcommand{\dar}{\downarrow}
\newcommand{\down}{\downarrow}

\newcommand{\ocM}{\overline{{\mathcal M}}}

\newcommand{\smooth}{{\operatorname{sm}}}
\newcommand{\sing}{_{\operatorname{sing}}}
\newcommand{\bal}{{\operatorname{bal}}}

\newcommand{\red}{{\operatorname{red}}}
\newcommand{\norm}{^{\operatorname{norm}}}
\newcommand{\double}{\genfrac..{0pt}1
{\raise -1pt\hbox{$\scriptstyle\longrightarrow$}}{\raise 3pt\hbox
{$\scriptstyle\longrightarrow$}}} 
\newcommand{\setmin}{\,\protect%
\begin{picture}(8,3.5)\qbezier(1,3.5)(4,2.)(7,.5)\end{picture}\,}
\newcommand{\curveuparrow}{%
\begin{picture}(4,10)%
\put(0,10){\vector(-1,2){0}}%
\qbezier(0,0)(4,5)(0,10)
\end{picture}}

\renewcommand{\setminus}{\setmin}

\newcommand{\mini}{{%
\begin{picture}(4,8)%
\put(2,4){\circle*{3}}%
\end{picture}}}

\newcommand{\KO}[4]{{\cK_{#1,#2}(#3,#4)}}
\newcommand{\KOB}[4]{{\cK^{\bal}_{#1,#2}(#3,#4)}}
\newcommand{\ko}[4]{{\bK_{#1,#2}(#3,#4)}}

\newcommand{\TSM}{\KO{g}{n}{\cM}{d}}
\newcommand{\tsm}{\ko{g}{n}{\cM}{d}}
\newcommand{\SM}{\KO{g}{n}{\bM}{d}}
\newcommand{\sm}{\ko{g}{n}{\bM}{d}}
\newcommand{\TSMB}{\KOB{g}{n}{\cM}{d}}

\def\sp{_{\rm sp}}
\def\eqdef{\mathrel{\mathop=\limits^{\rm def}}}

\def \fcoh{\mathop{\rm F\,Coh}\nolimits}
\def \coh{\mathop{\rm Coh}\nolimits}
\def\xet{{\cX}_{\rm \acute et}}

\def\cech{\v Cech}

\def\cechh{{\check{\rm H}}}

\def\cechc{{\check{\rm C}}}
\def\et{_{\rm \acute et}}

\def\mmu{{\protect\boldsymbol{\mu}}}

\def\gen{_{\rm gen}}
\def\rest{|_}

\def\tototi{\mathbin{\mathop{\otimes}\limits^{\raise-1pt\hbox
{$\scriptscriptstyle {\rm L}$}}}}

\def\indlim{\mathop{\vrule width0pt height7pt depth
4pt\smash{\lim\limits_{\raise 1pt\hbox to 14.5pt
{\rightarrowfill}}}}}
\def\projlim{\mathop{\vrule width0pt height7pt depth
4pt\smash{\lim\limits_{\raise 1pt\hbox to 14.5pt
{\leftarrowfill}}}}}

\begin{document}
\title[Stable maps]{Compactifying the space of stable maps}
\author[D. Abramovich]{Dan Abramovich}
\thanks{D.A. Partially supported by NSF grant DMS-9700520 and by an Alfred
	P. Sloan research fellowship}  
\address{Department of Mathematics\\ Boston University\\ 111 Cummington
	 Street\\ Boston, MA 02215\\ USA} 
\email{abrmovic@math.bu.edu}
\author[A. Vistoli]{Angelo Vistoli}
\thanks{A.V.  partially supported by the University of
	Bologna, funds for selected research topics.}
\address{Dipartimento di Matematica\\ Universit\`a di Bologna\\Piazza di Porta
	 San Donato 5\\ 40127 Bologna\\ Italy}
\email{vistoli@dm.unibo.it}
\date{\thisdate}

\begin{abstract} 
In this paper we study a  notion of {\em twisted stable map},
from a curve to a tame Deligne--Mumford stack, which generalizes the
well-known notion  of stable map to a projective variety.
\end{abstract}

\maketitle

\setcounter{tocdepth}{1}
\tableofcontents
\section{Introduction}

We fix a noetherian base scheme $\bbS$.

\subsection{The problem of moduli of families} Consider a Deligne-Mumford stack
$\cM$ 
admitting a projective coarse moduli scheme $\bM\subset \PP^N$. Given a curve
$C$, it is often natural to consider morphisms $f:C \to \cM$ (or equivalently,
objects $f\in \cM(C)$): in case $\cM$ is the moduli of geometric objects, these
morphisms correspond to families over $C$. For example, if
$\cM=\ocM_{\gamma}$, then 
morphisms $f:C\to \cM$ correspond to families of stable curves of genus
$\gamma$ over $C$; and if 
$\cM = \cB G$ we get principal $G$-bundles over $C$. It is interesting to study
moduli of such objects; moreover, it is natural 
to study such moduli as $C$ varies, and find a natural compactification for
such moduli. 

One approach is suggested by Kontsevich's moduli of stable maps.

\subsection{Stable maps}
 First consider a projective scheme $\bM \to \bbS$ with a fixed ample sheaf
 $\cO_\bM(1)$.
 Given integers
$g,n,d$, it is known that 
there exists a proper algebraic stack  $\SM$ of
stable, $n$-pointed maps of genus $g$ and degree $d$ into $\bM$ (see
\cite{Kontsevich}, \cite{Behrend-Manin}, \cite{Fulton-Pandharipande},
\cite{A-Oort}, where the notation $\overline\cM_{g,n}(\bM,d)$ is used). This 
stack admits a {\em projective} coarse moduli space $\sm$. If
one avoids ``small'' residue characteristics in $\bbS$, which depend
on $g,n,d$ and $\bM$, then this
stack is in fact a proper Deligne--Mumford stack.

\subsection{Stable maps into stacks}
 Now fix a proper Deligne--Mumford
stack $\cM \to \bbS$ admitting a projective coarse moduli space $\bM
\to \bbS$ on which we fix an ample sheaf as above. We further assume that $\cM$
is {\em 
tame}, that is, for  
any geometric point $s\colon  \Spec \Omega \to \cM$, the group 
$\Aut_{\Spec \Omega}(s)$ has order prime to the characteristic of the
algebraically closed field $\Omega$.

It is tempting to mimic Kontsevich's construction as
follows: let $C$ be a nodal projective connected curve; then a morphism $C \to
\cM$  is said to be a {\em stable map} of degree $d$ if the associated morphism
to the coarse moduli scheme $C \to \bM$ is a stable map of degree $d$.  

It follows from our results below that the category of stable maps into $\cM$
is  a 
Deligne-Mumford stack. A somewhat surprising point is that it is not
complete.

To see this, we fix $g=2$ and consider the specific case of $\cB G$ with $G=
(\ZZ/3\ZZ)^4$. Any 
smooth curve of genus $2$ admits a {connected} principal $G$ bundle,
corresponding of a surjection $H_1(C, \ZZ) \to G$. If we let $C$ degenerate to
a nodal curve $C_0$ of geometric genus $1$, then   $H_1(C_0, \ZZ) \simeq
\ZZ^3$, and since there is no surjection $\ZZ^3 \to G$, there is no connected
principal $G$-bundle over $C_0$. This means that there is no limiting stable
map  $C_0 \to \cB G$.

\subsection{Twisted stable maps}

Our main goal here is to correct this deficiency. In
order to do  so, we will enlarge the category of stable maps
into $\cM$. The source curve
${\cC}$ of a new stable map ${\cC} \to \cM$ will
acquire an orbispace structure at its nodes. Specifically, we
endow it  with the structure of a Deligne-Mumford stack.

It is not hard to see how these orbispace structures come
about. Let
$S$ be the spectrum of a discrete valuation ring $R$ of pure characteristic 0,
with quotient field $K$, and let
$C_K\to \eta\in S$ be a nodal curve over the generic point, together
with a map $C_K \to \cM$ of degree $d$, whose associated map $C_K \to \bM$ is
stable. We can 
exploit the fact that
${\cK}_{g,0}(\bM,d)$ is complete; after a ramified base
change on
$S$ the induced map $C_K \to
\bM$ will extend to a stable map $C \to \bM$ over $S$. Let
$C_\smooth$ be the smooth locus of the morphism $C \to S$;
Abhyankar's lemma, plus a fundamental purity lemma (see
\ref{Lem:purity-lemma} below) shows that after a suitable base change we can
extend the map $C_K \to
\cM$ to a map $C_\smooth \to \cM$; in fact the purity lemma fails
to apply only at the ``new'' nodes of the central fiber, namely those which are
not in the closure of nodes in the generic fiber. On the other hand, if
$p\in C$ is such a node, then on an
\'etale neighborhood $U$ of
$p$, the curve $C$ looks like $$uv = t^r,$$ where $t$ is the
parameter on the base. By taking
$r$-th roots:
$$u = u_1^r;\, v = v_1^r$$
we have a {\em nonsingular} cover
$V_0\to U$ where $V_0$ is defined by
$u_1v_1 = t$. The purity lemma applies to $V_0$,
so the composition ${V_0}_K \to C_K \to \cM$ extends over all
of $V_0$. There is a minimal intermediate cover $V_0\to
V\to U$ such that the family extends already over
$V$; this $V$ will be of the form $xy = t^{r/m}$, and the
map $V \to U$ is given by $u = x^m$,
$v = y^m$. Furthermore, there is an action of the group
$\bmu_m$ of roots of 1, under which $\alpha\in \bmu_m$ sends $x$
to $\alpha x$ and $y$ to $\alpha^{-1} y$, and
$V/ \bmu_m = U$. This gives the orbispace structure ${\cC}$
over $C$, and the map $C_K \to
\cM$ extends to a map ${\cC} \to \cM$.

This gives the flavor of our definition.

We define a
category $\TSM$, fibered over $\cS ch/\bbS$, of {\em twisted
stable $n$-pointed maps  $\cC 
\to \cM$ of genus $g$ and degree $d$}. This category is given in two
equivalent realizations: one as a 
category of stable {\em twisted $\cM$-valued objects} over  nodal
pointed curves endowed with atlases of orbispace charts (see
definition \ref{Def:stable-twisted-object}); the other as
a category of representable maps from pointed nodal Deligne--Mumford stacks
into $\cM$, such that the map on coarse moduli spaces is stable (see
definition \ref{Def:twisted-stable-map}). Both
realizations are used in proving our main theorem:

\begin{theorem}\label{Th:stable-maps}
\begin{enumerate}
\item 
The category $\TSM$ is a proper algebraic stack. 
\item The coarse moduli space
$\tsm$ of $\TSM$ is projective. 
\item There is a commutative diagram 
$$\begin{array}{ccc} \TSM & \to & \SM \\ 
\dar & & \dar \\ \tsm & \to & \sm
\end{array}
$$
where the top arrow is proper, quasifinite, relatively of
Deligne--Mumford type and tame, and the bottom arrow is finite. In particular,
if $\SM$ is a Deligne--Mumford stack, then so is
$\TSM$.
\end{enumerate} 
\end{theorem}

\subsection{Some applications and directions of further work} In our paper
\cite{A-V:fibered-surfaces} we 
studied the situation where $\cM = \overline\cM_{\gamma,\nu}$, which gives a
complete 
moduli for {\em fibered surfaces}. Further applications which we hope to
discuss in future papers  include: (1) The case where $\cM$ is the classifying
space of a finite group allows one to improve on the spaces of admssible covers
and to give moduli compactifications of Mumford's spaces of curves with level
structures. This is indicated in our \cite{A-V:stable-maps}, and is the subject
of investigation joint with Alessio Corti and Johan de Jong. This approach is
closely related to the work of Wewers \cite{Wewers}. (2) A similar
reasoning applies to  curves with  $r$-spin structures, e.g. theta
characteristics (joint investigation with Tyler Jarvis).  (3) The recursive
nature of the  
theorem allows one to construct both minimal models and stable reduction for
pluri-fibered varieties. This is related to recent work of Mochizuki
\cite{Mochizuki}. 

In this paper we verify that $\TSM$ is a proper stack by going through the
conditions one by one. It may be worthwhile to  develop a theory
of Grothendieck Quot-stacks and deduce our results from such a theory. It seems
likely that some of our methods could be useful for developing such a theory.

We were told by Maxim Kontsevich that he has also discovered the stack of
twisted stable maps, but has not written down the theory. His motivation was in
the direction of Gromov--Witten invariants of stacks.

\subsection{Acknowledgements} We would like to thank Kai Behrend, Larry Breen,
Barbara Fantechi, Ofer Gabber, 
Johan de Jong,  Maxim Kontsevich, and Rahul Pandharipande,  for helpful
discussions.  We are grateful to Laurent Moret-Bailly for providing us
with a preprint of the book \cite{L-MB} before it appeared. The first author
thanks the Max Planck Institute f\"ur Mathematik in Bonn for a  
visiting period which helped in putting this paper together.

\section{Generalities on stacks}
\subsection{Criteria for a Deligne--Mumford stack}\label{Sec:stack-criteria}
We refer the reader to \cite{Artin} and \cite{L-MB} for a general
discussion of algebraic stacks (generalizing \cite{Deligne-Mumford}), and to
the appendix in 
\cite{Vistoli:chow-stack} 
for an introduction. We spell out the conditions here, as we will follow them
closely in the paper. We are given a
category $\cX$ along with a functor $\cX \to \cS ch/\bbS$. We assume

\begin{enumerate}
\item \label{It:fibered-by-groupoids}
 $\cX \to \cS ch/\bbS$ is {\em fibered in groupoids} (see  \cite{Artin}, 
\S 1, (a) and (b) or \cite{L-MB}, Definition 2.1). This means:
\begin{enumerate} \item \label{It:pullbacks-exist}
 for any morphism of schemes $T \to T'$ and any object
$\xi'\in \cX(T')$ there is an object $\xi\in \cX(T)$ and an arrow $\xi \to
\xi'$ over $T\to T'$; and
\item \label{It:pullbacks-unique} For any diagram of schemes 
$$\begin{array}{ccccc}T_1 & &\lrar&  &  T_2 \\
			&\searrow&&\swarrow& \\
			&& T_3 ,&&
\end{array}$$ and any objects $\xi_i \in \cX(T_i)$ sitting in a compatible
diagram
$$\begin{array}{ccccc}\xi_1 & &&  &  \xi_2 \\
			&\searrow&&\swarrow& \\
			&& \xi_3 ,&&
\end{array}$$
there is a unique arrow $\xi_1 \to  \xi_2$ over $T_1 \to  T_2$ making the
diagram commutative.
\end{enumerate}

We remark that this condition is automatic for moduli problems, where
$\cX$ is a category of families with morphisms given by fiber diagrams.

\item \label{It:stack} $\cX \to \cS ch/\bbS$ is a {\em stack}, namely
 \begin{enumerate}
	\item \label{It:Isom-sheaf} the $\Isom$ functors are sheaves in
		the \'etale topology, and 
	\item \label{It:etale-descent} any \'etale descent datum
	for obejects of $\cX$ is effective.
 \end{enumerate}
 See \cite{Artin}, 1.1 or \cite{L-MB}, Definition 3.1.

\item The stack $\cX \to \cS ch/\bbS$ is {\em algebraic}, namely:

 \begin{enumerate}
	\item \label{It:Isom-representable} the $\Isom$ functors are
	representable by algebraic spaces locally of finite type, and
	\item \label{It:parametrization} There is a scheme $X$,
	locally of finite type, and  a 
	smooth and surjective morphism $X \to \cX$.
 \end{enumerate}
See \cite{Artin}, Definition 5.1. This differs slightly from
	\cite{L-MB}, where one assumes in addition that $\cX \to
	\cX\times\cX$ is separated. 
Notice that (\ref{It:Isom-representable}) implies (\ref{It:Isom-sheaf}). 

These last two conditions are often the most difficult to verify.
For the last one, M. Artin has devised a set of criteria for constructing
$X\to \cX$ by algebraization of  formal deformation spaces (see
\cite{Artin}, Corollary 5.2). Thus, {\em in case $\bbS$ is of finite type over
a field or an excellent Dedekind domain,} condition
\ref{It:parametrization} holds if 
\begin{enumerate}[(A)]
\item \label{It:limit-preserving}
$\cX$ is limit preserving (see \cite{Artin}, \S 1);
\item \label{It:formal-compatible}
$\cX$ is compatible with formal completions (see \cite{Artin},
5.2 (3)); 
\item \label{It:pro-rep}
Schlessinger's conditions for pro-representability of the
deformation functors hold (see \cite{Artin}, (2.2) and (2.5)); and
\item \label{It:def-obs}
There exists an obstruction theory for $\cX$ (see \cite{Artin},
(2.6)) such that 
\begin{enumerate}
 \item  \label{It:def-obs-etale}
 the deformation and obstruction theory is compatible with
 	\'etale localization (\cite{Artin}, 4.1 (i));
 \item \label{It:def-obs-formal}
	 the deformation theory is compatible with formal completions
 	(\cite{Artin}, 4.1 (ii)); and 
 \item  \label{It:def-obs-constructible}
	 the deformation and obstruction theory is constructible
 	(\cite{Artin}, 4.1 (iii)). 
\end{enumerate}

\end{enumerate}

\end{enumerate}
Furthermore, we say that $\cX$ is a {\em Deligne--Mumford stack} if we can
choose $X \to \cX$ as in (\ref{It:parametrization}) to be {\em \'etale}. This
holds if and only if the diagonal $\cX \to \cX \times \cX$ 
is unramified.  A morphism $\cX \to \cX_1$ is {\em of Deligne--Mumford type}
if for any scheme $Y$ and morphism $Y \to \cX_1$ the stack $Y
\times_{\cX_1} \cX$ is a Deligne--Mumford stack.  

For the notion of {\em properness} of an algebraic stack see \cite{L-MB},
Chapter 7. Thus a stack $\cX\to \bbS$ is proper if it is separated, of finite
type and universally closed. In  \cite{L-MB} it is noted that the weak
valuative criterion for properness using traits might be insufficient for 
properness. However, in case $\cX$ has {\em finite diagonal}, it is shown in
\cite{Edidin}, Corollary 4.1, that there exists  a finite  
surjective morphism  from  a scheme $Y \to \cX$. In such a case the usual weak
valuative criterion suffices.

\subsection{Coarse moduli spaces}

Recall the following result:

\begin{theorem}[Keel-Mori \cite{Keel-Mori}]\label{Th:moduli-space} Let $\cX$ be
an algebraic stack 
with finite diagonal over a scheme $S$. There exists an algebraic space $\bX$
and a  morphism $\cX \to \bX$ such that
\begin{enumerate}
\item\label{It:moduli-space-finite} $\cX \to \bX$ is proper and quasifinite; 
\item\label{It:moduli-space-points} if $k$ is an algebraically closed field,
then $\cX(k)/\Isom \to 
\bX(k)$ is a bijection.
\item\label{It:moduli-space-categorical}  whenever $Y\to S$ is an algebraic
space and $\cX 
\to Y$ is a morphism, then the morphism factors uniquely as $\cX \to \bX \to
Y$; more generally
\item\label{It:moduli-space-universal} whenever $S' \to S$ is a flat morphism
of schemes, and  whenever $Y\to 
S'$ is  
an algebraic space and $\cX\times _S S' 
\to Y$ is a morphism, then the morphism factors uniquely as $\cX\times _S S'
\to \bX\times _S S' \to  
Y$; in particular
\item\label{It:moduli-space-sheaf} 
$\pi_* \cO_{\cX} = \cO_{\bX}$
\end{enumerate}
\end{theorem}

Recall that an algebraic space $\bX$ along with a morphism $\pi\colon \cX
\to \bX$ satisfying properties \ref{It:moduli-space-points} and
\ref{It:moduli-space-categorical} is called a {\em coarse moduli space} (or
just {\em moduli space}). In particular, the theorem of Keel and Mori shows
that coarse moduli spaces of algebraic stacks with finite diagonal exist. 
Moreover, from \ref{It:moduli-space-universal} and \ref{It:moduli-space-sheaf}
above we have that the formation of a coarse moduli space behaves well under
flat base change: 

\begin{lemma}\label{Lem:charac-cms} Let ${\cX} \to X$ be a proper
quasifinite morphism, where
${\cX}$ is a Deligne--Mumford stack and $X$ a
noetherian scheme. Let
$X'
\to X$ be a flat morphism of schemes, and denote ${\cX}' = X'
\times_X {\cX}$. 
\begin{enumerate}
\item  If $X$ is the moduli space of ${\cX}$, then $X'$
is the moduli space of ${\cX}'$.

\item  If $X' \to X$ is also surjective and $X'$ is the
moduli space of ${\cX}'$, then $X$ is the moduli space of
${\cX}$.
\end{enumerate}
\end{lemma}

\proof Given a proper quasifinite morphism $\pi\colon  {\cX} \to
X$, then it exhibits $X$ as a moduli space if and only if
$\pi_*{\cO}_{\cX} = {\cO}_X$. If $R\ \double  \
U$ is an \'etale presentation of ${\cX}$, $f \colon   U \to
X$ and $g \colon   R \to X$ the induced morphisms, then this
condition is equivalent to the exactness of the sequence
$$
0 \lrar {\cO}_X \lrar f_* {\cO}_U \double g_* {\cO}_R.
$$
From this the statement follows. \endproof

The prototypical example of a  moduli space is given by a group quotient: let
$V$ be a 
scheme and $\Gamma$ a finite group acting on $V$. The morphism
$[V/\Gamma] \to V/\Gamma$ exhibits the quotient space $V/\Gamma$ as the moduli
space of the 
stack $[V/\Gamma]$. The following well-known lemma shows that \'etale-locally,
the moduli space of any Deligne--Mumford stack is of this form.

\begin{lemma}\label{Lem:locally-quotient} Let $\cX$ be a separated
Deligne--Mumford stack, and $X$ its  
coarse moduli space. There is
an \'etale covering $\{X_\alpha \to X\}$, such that for each
$\alpha$ there is a scheme $U_\alpha$ and a finite group
$\Gamma_\alpha$ acting on $U_\alpha$, with the property
that the pullback
${\cX}\times_XX _\alpha$ is isomorphic to the
stack-theoretic quotient $[U_\alpha/\Gamma_\alpha]$.
\end{lemma}

{\em Sketch of proof.}
 Let $x_0$ be a geometric point of $X$. Denote by $X^\sh$ the
spectrum of the strict henselization of $X$ at the point
$x_0$, and ${\cX}^\sh = {\cX}\times_X X^\sh$. If $V \to {\cX}$ is an
\'etale morphism, with
$V$ a scheme, having $x_0$ in its image, there is a component
$U$ of the pullback $V \times_X X^\sh$ which is finite over
$X^\sh$. Denote $R = U \times_{{\cX}^\sh}U$. We have that
under the first projection $R \to U$, the scheme $R$ splits as a disjoint 
union of copies of $U$. Let $\Gamma$ be the set of
connected components of $R$, so that $R$ is isomorphic to $U
\times\Gamma$. Then the product $R \times_U R \to R$ induces a
group structure on $\Gamma$, and the second projection $R
\simeq U \times\Gamma \to U$ defines a group action of\/
$\Gamma$ on $U$, such that ${\cX}^\sh$ is the quotient
$U/\Gamma$.

The statement follows from standard limit arguments.\endproof

\subsection{Tame stacks and their coarse moduli spaces}
\begin{definition}\begin{enumerate} 
\item A Deligne--Mumford  stack $\cX$ is said to be {\em tame} if 
for
any geometric point $s\colon  \Spec \Omega \to \cX$, the group 
$\Aut_{\Spec \Omega}(s)$ has order prime to the characteristic of the
algebraically closed field $\Omega$. 
\item A morphism $\cX \to \cX_1$ of algebraic stacks is said to be {\em tame}
if for any scheme $Y$ and morphism $Y \to \cX_1$ the stack $Y 
\times_{\cX_1} \cX$ is a tame Deligne--Mumford stack.  
\end{enumerate}
\end{definition}

A closely related notion is the following:

\begin{definition}
 An action of a finite group $\Gamma$ on a scheme $V$ is said to be {\em
tame} if  for
any geometric point $s\colon  \Spec \Omega \to V$, the group $\Stab (s)$ has
order 
prime to the characteristic of $\Omega$.
\end{definition}

The reader can verify that a separated Deligne--Mumford   stack is tame if and
only if the actions of the groups $\Gamma_\alpha$ on $V_\alpha$ in the previous
lemma are tame.

In case $\cX$ is tame, the formation of coarse moduli spaces
commutes with arbitrary morphisms: 

\begin{lemma}\label{Lem:tame-cms-pullback} Let ${\cX}$ be a tame
Deligne--Mumford stack, ${\cX} \to X$ its moduli
space. If $X' \to X$ is any morphism of schemes, then $X'$ is
the moduli space of the fiber product $X'
\times_X {\cX}$. Moreover, if $X'$ is reduced, then it is also the
moduli space of $(X'
\times_X {\cX})_\red$.
\end{lemma}

\proof  By Lemma \ref{Lem:charac-cms}, this is a local condition in
the \'etale topology of $X$, so we may assume that ${\cX}$
is a quotient stack of type $[Y / \Gamma]$, where $\Gamma$ is a
finite group acting on an affine scheme $Y = \Spec R$. Moreover, since $\cX$ is
tame, we may assume that the order of $\Gamma$ is prime to all residue
characteristics. Then 
$X= \Spec R^ \Gamma$; if $X' = \Spec S$, then the statement
is equivalent to the map $S \to \left(R \otimes_{R^ \Gamma}
S\right)^ \Gamma$ being an isomorphism. This (well known) fact can be shown as
follows: Recall that for any $R^\Gamma[\Gamma]$-module $M$ the homomorphism 
\begin{eqnarray*}  
M &\stackrel{\Av_M}{\lrar}& M^ \Gamma \\
   y & \mapsto & \frac{1}{|\Gamma|} \sum_{\gamma\in \Gamma} \gamma \cdot y
\end{eqnarray*} 
is a projector exhibiting $M^\Gamma$ as a direct summand in
$M$. Thus the induced morphism $$\Av_R\otimes \Id_S\colon  R\otimes_{R^ \Gamma}
S \to 
R^\Gamma\otimes_{R^\Gamma}S = S$$ shows that $S \to \left(R \otimes_{R^ \Gamma}
S\right)^ \Gamma$ is injective. The morphism $\Av_R\otimes \Id_S$ is a lifting
of $$\Av_{R 
\otimes_{R^ \Gamma} S}\colon  R \otimes_{R^ \Gamma} S \to \left(R \otimes_{R^
\Gamma} 
S\right)^ \Gamma, $$ which is surjective. 

This shows that $X'$ is
the moduli space of the fiber product $X'
\times_X {\cX}$.
The statement about $(X'
\times_X {\cX})_\red$ is immediate. This proves the result.
\endproof

Let $\cX$ be a separated {\em tame} stack with coarse moduli scheme
$\bX$. Consider  
the projection $\pi\colon  \cX \to \bX$. The functor
$\pi_*$ carries sheaves of ${\cO}_{\cX}$-modules to
sheaves of ${\cO}_\bX$-modules.

\begin{lemma}\label{Lem:cms-exact} The functor $\pi_*$ carries
quasicoherent sheaves to quasicoherent sheaves, coherent
sheaves to coherent sheaves, and is exact.
\end{lemma}

\proof The question is local in the \'etale topology on $\bX$,
so we may  assume that ${\cX}$ is of
the form $[V/\Gamma]$, where $V$ is a scheme and $\Gamma$ a
finite group of order prime to all resudue characteristics, in particular $\bX
= V/\Gamma$. Now sheaves on ${\cX}$ correspond to 
equivariant sheaves on $V$. Denote by $q\colon V \to \bX$ the projection. If
$\cE$ is 
a sheaf on $\cX$ corresponding to a $\Gamma$-equivariant sheaf $\tilde\cE$ on
$V$, then  $\pi_*\cE = (q_*\tilde\cE)^\Gamma$, which, by the tameness
assumption,  is a direct summand in $q_*\tilde\cE$. From this the statement
follows.\endproof

\subsection{Purity lemma}

We recall the following {\em purity lemma} from \cite{A-V:fibered-surfaces}:

\begin{lemma}\label{Lem:purity-lemma}
 Let $\cM$ be a
separated 
Deligne - Mumford 
stack, $\cM\to
\bM$ the coarse moduli space. Let $X$ be a separated scheme of
dimension 2 satisfying Serre's condition $S_2$. Let $P\subset X$ be a finite
subset consisting of closed points, $U=X\setmin P$. Assume that the local
fundamental groups of $U$
around the points of $P$ are trivial.

Let $f\colon  X \to \bM$ be a morphism. Suppose there is a lifting
$\tilde{f}_U\colon U \to \cM$:
\begin{equation}\begin{diagram} \node[3]{\cM}\arrow{s} \\
\node{U}\arrow{ene,t}{\tilde{f}_U}\arrow{e} 
\node{X}\arrow{e,t}{f} \node{\bM} \end{diagram}
\end{equation}

Then the lifting extends  to $X$:
\begin{equation}\begin{diagram} \node[3]{\cM}\arrow{s} \\
\node{U}\arrow{ene,t}{\tilde{f}_U}\arrow{e} 
\node{X}\arrow{e,t}{f}\arrow{ne,t,1,..}{\tilde{f}} \node{\bM}
\end{diagram}
\end{equation}
and $\tilde{f}$ is unique up to a unique isomorphism.
\end{lemma}

\proof
By the descent axiom for $\cM$ (see
\ref{Sec:stack-criteria}-(\ref{It:stack})) the problem is local in the
\'etale topology, so we may 
replace $X$ and $\bM$ with the spectra of their strict henselizations
at a geometric point; then we can also assume that we have a universal
deformation space
$V\to \cM$ which is {\em finite}. Now $U$ is the complement of the
closed point, $U$ maps to
$\cM$, and the pullback of $V$ to $U$ is finite and \'etale, so it has
a section, because $U$ is simply connected; consider the corresponding map
$U\to V$. Let $Y$ be the scheme-theoretic closure of the graph of this map
in $X\times_\bM V$. Then $Y\to X$ is finite and is an isomorphism on $U$.
Since $X$ satisfies  $S_2$, the morphism $Y\to X$ is an isomorphism. \qed

\begin{remark}
The reader can verify that the statement and proof work in higher
dimension. See also  related lemmas in \cite{Mochizuki}.
\end{remark}

\begin{corollary}
Let $X$ be a smooth surface over a field, $p\in X$ a closed point with
complement $U$. 
Let $X\to \bM$ and $\, U\to \cM$ be as in the purity lemma. Then
there is a lifting $X \to \cM$.
\end{corollary}

\begin{corollary}
Let $X$ be a normal crossings surface over a field $k$, namely a surface which
is  \'etale locally isomorphic to $\Spec k[u,v,t]/(uv)$. Let $p\in X$ a closed
point with complement $U$. Let $X\to \bM$ and $\, U\to \cM$ be as
in 
the purity lemma. Then there is a lifting $X \to \cM$.
\end{corollary}

\proof
In both cases $X$ satisfies condition $S_2$ and the local fundamental group
around $p$ is trivial, hence the purity lemma applies. \endproof

\subsection{Descent of equivariant objects}
\begin{lemma}\label{Lem:descent} Let $R$ be a local ring with residue field
$k$, let $U = \Spec 
R$, $u_0 = \Spec k$,  and let $\eta$ be
an object of ${\cM}(U)$. Assume we have a pair of compatible actions of a
finite group $\Gamma$ on 
$R$ and on $\eta$, in such a way that the induced actions
of\/ $\Gamma$ on $k$ and on the pullback $\eta_0=\eta|_{u_0}$ are trivial. Then
there exists an  
object $\eta'$ of $\cM$ on the quotient $U/ \Gamma= \Spec(R^ \Gamma)$,   
and a $\Gamma$-invariant lifting $\eta \to\eta'$ of the
projection $U \to U/
\Gamma$. Furthermore, if $\eta''$ is another such object
over $U/ \Gamma$, there is a unique isomorphism $\eta' \simeq
\eta''$ over the identity of $U/ \Gamma$, which is compatible
with the two arrows $\eta \to\eta'$ and $\eta\to\eta''$.
\end{lemma}

As a consequence of the unicity statement, suppose that we
have a triple $( \alpha, \beta, \gamma)$, where $\gamma\colon 
\Gamma\simeq\Gamma$ is a group isomorphism, and $\alpha\colon 
\eta\simeq\eta$ and $\beta\colon  U \simeq U$ are compatible
$\gamma$-equivariant isomorphisms. Then the given arrow
$\eta\to\eta'$ and its composition with $\alpha$ both satisfy
the conditions of the lemma, so there is an induced
isomorphism $\overline \alpha\colon  \eta' \simeq\eta'$.

\begin{corollary}\label{Cor:descent-to-quotient} Let $R,U,k,u_0,\eta,\eta_0$ be
as in the previous lemma. 
Let 
$G$ be a finite group acting
compatibly on $R$ and on $\eta$. Let
$\Gamma$ the normal subgroup of $G$ consisting of elements
acting on $k$ and $\eta_0$ as the identity. Then there is a $G/
\Gamma$-equivariant object $\eta'$ on the quotient $U/
\Gamma$, and a $G$ equivariant arrow $\eta\to\eta'$
compatible with the projection $U \to U/ \Gamma$.
\end{corollary}

{\em Proof of the corollary.} The action is defined as follows. If $g$ is an 
element 
of $G$, call $\alpha\colon  \eta\simeq\eta$ and $\beta\colon  U
\simeq U$ the induced arrows, and $\gamma\colon 
\Gamma\simeq\Gamma$ the conjugation by $g$. Then the image of
$g$ in $G/ \Gamma$ acts on $\eta'$ via the isomorphism
$\overline \alpha\colon  \eta' \simeq\eta''$ defined above.
One checks easily that this defined an action with the
required properties.\endproof

{\em Proof of the Lemma.} First note that if $R^\sh$ is
the strict henselization of $R$, the condition on the action
of\/ $\Gamma$ allows to lift it to $R^\sh$. Also, the statement
that we are trying to prove is local in the \'etale topology,
so by standard limit arguments we can assume that $R$ is
strictly henselian. Replacing $\bM$ by the
spectrum of the strict henselization of its local ring at the
image of the closed point of $R$, we can assume that ${\cM}$ is of the
form $[V/H]$, where $V$ is a scheme and $H$ is a   
finite group. Then the object $\eta$ corresponds to a
principal $H$-bundle $P \to U$, on which
$\Gamma$ acts compatibly with the action of\/ $\Gamma$ on $U$,
and an
$H$-equivariant and $\Gamma$-invariant morphism $P \to V$.
Since $U$ is strictly henselian, the bundle $P \to U$ is trivial, so $P$ is a
disjoint union 
of copies of $\Spec R$, and the group
$\Gamma$ permutes these copies;
furthermore the hypothesis on the action of\/ $\Gamma$ on
the closed fiber over the residue field insures that $\Gamma$
sends each component into itself. The thesis follows
easily.\endproof

\section{Twisted objects}

\subsection{Divisorially marked curves}

The following definition is a local version of the standard
definition of pointed curve; its advantage is that it is stable
under localization in the \'etale topology.

\begin{definition} A {\em divisorially $n$-marked
nodal curve}, or simply {\em $n$-marked curve} $(U \to S, \Sigma_i)$,
consists of a nodal 
curve $\pi\colon  U
\to S$, together with a sequence of $n$ pairwise disjoint
closed subschemes
$\Sigma_1, \ldots, \Sigma_n\subset U$  whose
supports do not contain any of the singular points of the
fibers of
$\pi$, and such that the projections $\Sigma_i \to S$ are
\'etale. (Any of the subschemes $\Sigma_i$ may be empty.)
\end{definition}

If more than one curve  is considered, we will often use the notation
$\Sigma_i^U$ to specify the curve $U$. On the other hand, 
we will often omit  the subschemes $\Sigma_i^U$ from the notation $(U \to S,
\Sigma_i^U)$ if there is no risk of confusion.

A nodal $n$-pointed curve $C \to S$ is considered
an $n$-marked curve by taking as the
$\Sigma^C_i$ the images of the sections $S \to C$.

\begin{definition} If $(U \to S,\Sigma_i)$ is an $n$-marked nodal
curve, we define the {\it special locus\/} of $U$, denoted
 by $U\sp$, to be the union of the $\Sigma_i $
with the singular locus of the projection $U \to S$, with its natural scheme
 structure (this makes the 
projection $U\sp \to S$ unramified). The complement of
$U\sp$ will be called the {\it general locus\/} of $U$, and
denoted by $U\gen$.
\end{definition}

\begin{definition}If $(U \to S,\Sigma_i^U)$ is a marked curve, and $S'
\to S$ is an arbitrary morphism, we define the pullback to be $(U'\to
S',\Sigma^{U'}_i)$, where  $U' = S' \times_S 
U$,  and  $
\Sigma^{U'}_i = S' \times_S \Sigma^U_i$.
\end{definition}

\begin{definition} If $U \to S$ and $V \to S$ are
$n$-marked curves, a\/ {\em morphism of $n$-marked curves} $f
\colon   U \to V$ is a morphism of
$S$-schemes which sends each $\Sigma^U_i$ into $\Sigma^V_i$.

A morphism of $n$-marked curves $f \colon   U \to V$ is
called\/ {\em strict} if the support of
$f^{-1}(\Sigma^V_i)$ coincides with the support of
$\Sigma^U_i$ for all $i = 1, \ldots, n$, and similarly for the singular locus.
\end{definition}

One should notice that if a morphism of marked curves
$U \to V$ is strict, then there is an induced morphism of
curves $U\gen \to V\gen$. Furthermore, if $f
\colon   U \to V$ is strict and \'etale, then $f^{-1}(\Sigma^V_i)
= \Sigma^U_i$ scheme-theoretically.

\begin{definition} Let $(U \to S, \Sigma_i)$ be an $n$-marked curve,
$\Gamma$ a finite group. An {\em action} of\/
$\Gamma$ on $(U,\Sigma_i)$ is an action of\/ $\Gamma$ on $U$ as an
$S$-scheme,  
such that each element on $\Gamma$ acts via an automorphism
of $U$ as a marked curve on $S$.
\end{definition}

If $\Gamma$ is a finite group along with a {\em  tame} action on a marked curve
$U \to S$, then the quotient $U / \Gamma \to S$ can be given a
structure of marked curve by defining
$\Sigma^{U/ \Gamma}_i := \Sigma^U_i/ \Gamma \subseteq U/ \Gamma$.
The latter inclusion  holds because the orders of stabilizers in $\Gamma$ are
assumed to be prime to the
residue characteristics, so $\Sigma^U_i/ \Gamma$ is
indeed a subscheme of $U/\Gamma$.

Given a  morphism  $f \colon  U \to V$ of marked curves, and a tame action of a
 finite group $\Gamma$ on $U$, leaving $f$
invariant, then there is an induced morphism
$U/ \Gamma \to V$ of marked curves.


\begin{definition} Let $(U \to S, \Sigma_i)$ be an $n$-marked curve, with an 
action of a finite group 
$\Gamma$, and let $\cM$ be a Deligne--Mumford stack.
Given $\eta
\in \cM(U)$, an\/ {\em essential action} of\/ $\Gamma$
on $(\eta , U)$ is a pair of compatible actions of\/ $\Gamma$
on $\eta$ and on $(U\to S, \Sigma_i)$, with the property that if
$g$ is an element of\/ $\Gamma$ different from the
identity, and $u_0$ is a geometric point of $U$ fixed by
$g$, the automorphism of the pullback of $\eta$ to $u_0$
induced by $g$ is not trivial.
\end{definition}

\subsection{Generic objects and charts}

\begin{definition}
 Let  $C \to S$ be an
$n$-pointed nodal curve. A\/ {\em generic object} on $C$ is
an object of $\cM(C\gen)$.
\end{definition}

We will often write $(\xi, C)$ for a generic object $\xi$ on a curve $C$.

\begin{definition}
 Let $C \to S$ be an
$n$-pointed nodal curve, $\xi$ a generic object on $C$. A\/
{\em chart} $(U, \eta, \Gamma)$ for $\xi$ consists of the
following collection of data.
\begin{enumerate}
\item An $n$-marked curve $U \to S$, and
a strict morphism $\phi \colon  U \to C$.

\item An object $\eta$ of $\cM(U)$.

\item An arrow
$\eta|_{U\gen} \to \xi$ in $\cM$ 
 compatible with
the restriction $\phi|_{U\gen} \colon  U\gen \to C\gen$.

\item A finite group $\Gamma$.

\item An tame, essential action of\/ $\Gamma$ on $(\eta,U)$.

\end{enumerate}

Furthermore, we require that the following conditions be satisfied.
\begin{description}
\item[a] The actions of\/ $\Gamma$ leave the morphism $U \to
C$ and the arrow $\eta|_{U\gen} \to \xi$ invariant.

\item[b] The induced morphism $U/\Gamma \to C$ is \'etale.
\end{description}
\end{definition}

The following gives a local description of the action of $\Gamma$.

\begin{proposition} Let $(U, \eta, \Gamma)$ be a chart for
a generic object $\xi$ on a pointed nodal curve $C
\to S$. Then the action of\/ $\Gamma$ on $U\gen$ is free.

Furthermore, if $s_0$ is a geometric point of $S$ and $u_0$ a
nodal point of the fiber $U_{s_0}$ of $U$ over $s_0$, then

\begin{enumerate}
\item the stabilizer $\Gamma'$ of $u_0$ is a cyclic group
which sends each of the branches of $U_{s_0}$ to itself;

\item If $k$ is the order of\/ $\Gamma'$, then a generator
of\/ $\Gamma'$ acts on the tangent space of each branch by
multiplication with a primitive $k$-th root of 1.
\end{enumerate}

In particular, each nodal point of $U_{s_0}$ is sent to a
nodal point of $C_{s_0}$.
\end{proposition}

\proof The first statement follows from the definition of an
essential action, and the invariance of the arrow $\eta
\rest{U\gen} \to E$.

As for (1), observe that if the stabilizer
$\Gamma'$ of $u_0$ did not preserve the branches of $U_{s_0}$
then the quotient $U_{s_0}/\Gamma'$, which is \'etale at the
point $u_0$ over the fiber
$U_{s_0}$, would be smooth over $S$ at $u_0$, so $u_0$ would
be in the inverse image of
$U\gen$. From the first part the Proposition it would follow
that
$\Gamma'$ is trivial, a contradiction.

So $\Gamma'$ acts on each of the two branches individually.
The action on each branch must be faithful because it is free
on the complement of the set of nodes; this means that the
representation of\/ $\Gamma'$ in each of the tangent spaces to
the branches is faithful, and this implies the final
statement.\endproof

\begin{definition}\label{Def:balanced-chart}
 A chart is called\/ {\em balanced} if
for any nodal point of any geometric fiber of $U$, the two
roots of 1 describing the action of a generator of the
stabilizer on the tangent spaces to each branch of $U$ are
inverse to each other.
\end{definition}

\subsection{The transition scheme} \label{Sec:transition scheme}

Let $\xi$ be a generic object over a nodal curve $C\to S$ ,
$(U_1,\eta_1,\Gamma_1)$,
$(U_2,\eta_2,\Gamma_2)$ two charts; call ${\rm
pr}_i\colon  U_1\times_C U_2\to U_i$ the $i^{\rm th}$
projection. Consider the scheme
$$ I = {\mathop{\Isom}\limits_{U_1\times_C U_2}}({\rm
pr}_1^*\eta_1,{\rm pr}_2^*\eta_2)
$$
over $U_1\times_C U_2$ representing the functor of
isomorphisms of the two objects ${\rm pr}_1^*\eta_1$ and ${\rm
pr}_2^*\eta_2$. 

There is a section of $I$ over the inverse
image
$\widetilde U$ of
$C\gen$ in $U_1\times_C U_2$ which corresponds to the
isomorphism ${\rm pr}_1^*\eta_1|_{\widetilde U}\simeq {\rm
pr}_2^*\eta_2|_{\widetilde U}$ coming from the fact that
both
${\rm pr}_1^*\eta_1$ and ${\rm pr}_2^*\eta_2$ are pullbacks to
$\widetilde U$ of $\xi$. We will
call the scheme-theoretic closure $R$ of this section in $I$
the {\em transition  scheme\/} from $(U_1,\eta_1,\Gamma_1)$
to $(U_2,\eta_2,\Gamma_2)$; it comes equipped with two
projections $R\to U_1$ and
$R\to U_2$. 

There is also an action of
\/$\Gamma_1\times\Gamma_2$ on $I$, defined as follows. Let
$(\gamma_1,\gamma_2)\in \Gamma_1\times\Gamma_2$, and
$\phi\colon  {\rm pr}_1^*\eta_1\simeq {\rm pr}_2^*\eta_2$ an
isomorphism over $U_1\times_C U_2$; then define
$(\gamma_1,\gamma_2)\cdot \phi = \gamma_2\circ \phi\circ
\gamma_1^{-1}$. This action of\/
$\Gamma_1\times\Gamma_2$ on $I$ is compatible with the action
of\/ $\Gamma_1\times\Gamma_2$ on $U_1\times_C U_2$, and leaves
$R$ invariant. It follows from the definition of an essential
action that the action of\/ $\Gamma_1 = \Gamma_1\times \{1\}$
and $\Gamma_2 =
\{1\}\times\Gamma_2$ on $I$ is free.

\subsection{Compatibility of charts}
\begin{definition} Two charts $(U_1,\eta_1,\Gamma_1)$ and
$(U_2,\eta_2,\Gamma_2)$ are\/ {\em compatible} if their
transition scheme $R$ is \'etale over $U_1$ and $ U_2$.
\end{definition}

Let us analyze this definition. First of all, $R$ is
obviously \'etale over $(U_1)\gen$ and $(U_2)\gen$. Also,
since the maps $U_j \to C$ are strict, it is clear that the
inverse image of $\Sigma^{U_1}_i$ in $R$ is
set-theoretically equal to the inverse image of
$\Sigma^{U_2}_i$. If the two charts are compatible, this also
holds scheme-theoretically.

Now, start from two charts
$(U_1,\eta_1,\Gamma_1)$ and $(U_2,\eta_2,\Gamma_2)$.
Fix two geometric points
$$u_1\colon  \Spec \Omega
\to U_1 \quad\mbox{  and }\quad  u_2\colon  \Spec \Omega \to U_2$$ mapping to
the 
same geometric point $u_0\colon  \Spec
\Omega\to C$, and call $\Gamma'_j\subset \Gamma_j$ the
stabilizer of
$u_j$. Also call $ U_1^\sh$,
$ U_2^\sh$ and
$C^\sh$ the spectra of the strict henselizations of $ U_1 $, $
U_2 $ and $C$ at the points $u_1,u_2$ and $u_0$ respectively.
The action of\/ $\Gamma_j$ on $U_j$ induces an action of\/
$\Gamma'_j$ on $ U_j^\sh$. Also call $\eta_j^\sh$ the pullback of
$\eta_j$ to $ U_j^\sh$; there is an action of\/ $\Gamma'_j$ on
$\eta_j^\sh$ compatible with the action of\/ $\Gamma'_j$ on $
U_j$. The following essentially says that two chart are
compatible if for any choice of $u_1$ and $u_2$ the two
charts are locally isomorphic in the \'etale topology.

\begin{proposition} \label{Prp:compatibility-local}  The two charts are
compatible if and 
only if for any pair of geometric points $u_1$ and $u_2$ as
above there exist an isomorphism of groups $\theta\colon 
\Gamma'_1\simeq \Gamma'_2$, a $\theta$-equivariant
isomorphism $\phi\colon  U_1^\sh\simeq U_2^\sh$ of schemes over
$C^\sh$,  and a compatible $\theta$-equivariant isomorphism $\psi\colon
\eta_1^\sh\to 
\eta_2^\sh$.
\end{proposition}

\proof Consider the spectrum $(U_1\times_C U_2)^\sh$ of the
strict henselization of
$U_1\times_C U_2$ at the point $(u_1,u_2)\colon 
\Spec\Omega\to U_1\times_C U_2$, and call $R^\sh$ the pullback
of $R$ to $(U_1\times_C U_2)^\sh$. Assume that the two charts
are compatible. The action of\/ $\Gamma_1\times\Gamma_2$ on
$I$ described above induces an action of\/
$\Gamma'_1\times\Gamma'_2$ on $R^\sh$, compatible with the
action of\/ $\Gamma'_1\times\Gamma'_2$ on
$(U_1\times_C U_2)^\sh$. The action of\/ $\Gamma'_1 =
\Gamma'_1\times\{1\}$ on the inverse image of $C\gen$ in
$R^\sh$ is free, and its quotient is the inverse image of
$C\gen$ in $U_2^\sh$; but $R^\sh$ is finite and \'etale over
$U_2^\sh$, so the action of\/ $\Gamma'_1$ on all of $R^\sh$ is
free, and $R^\sh/\Gamma'_1 = U_2$. Analogously the action of\/
$\Gamma'_2$ on $R^\sh$ is free, and $R^\sh/\Gamma'_2 = U_1$.

Now, each of the connected components of $R^\sh$ maps
isomorphically onto both $U_1^\sh$ and
$U_2^\sh$, because $U_j^\sh$ is the spectrum of a strictly henselian
ring and the projection $R^\sh\to U_j^\sh$ is \'etale; this
implies in particular that the order of\/ $\Gamma_1$ is the
same as the number $k$ of connected components, and likewise
for
$\Gamma_2$. Fix one of these components, call it
$R_0^\sh$; then we get isomorphisms $R_0^\sh\simeq U_j^\sh$, which
yield an isomorphism $\phi\colon  U_1^\sh\simeq U_2^\sh$.

Call $\Gamma'$ the stabilizer of the component $R_0^\sh$
inside $\Gamma'_1\times\Gamma'_2$; the order of\/
$\Gamma'$ is at least
$|\Gamma'_1\times\Gamma'_2|/k = k^2/k = k$. But the action
of\/
$\Gamma'_2$ on $R^\sh$ is free, and so $\Gamma'\cap \Gamma_2 =
\{1\}$; this implies that the order of\/ $\Gamma'$ is $k$, and
the projection $\Gamma'\to \Gamma'_1$ is an isomorphism.
Likewise the projection $\Gamma'\to \Gamma'_2$ is an
isomorphism, so from these we get an isomorphism $\theta\colon 
\Gamma'_1\to
\Gamma'_2$, and it is easy to check that the isomorphism of
schemes
$\phi\colon  U_1^\sh\simeq U_2^\sh$ is $\theta$-equivariant.

There is also an isomorphism of the pullbacks of $\eta_1^\sh$
and $\eta_2^\sh$ to $R_0^\sh$, coming from the natural morphism
$R_0^\sh\to I$, which induces an isomorphism $\psi\colon 
\eta_1^\sh\to \eta_2^\sh$. This isomorphism is compatible with
$\phi$, and is it also $\theta$-equivariant.

Let us prove the converse. Suppose that there exist $\theta$,
$\phi$ and
$\psi$ as above. Then there is a morphism $\sigma\colon 
U_1^\sh\times\Gamma'_1\to I$ which sends a point
$(u_1,\gamma_1)$ of
$U_1^\sh\times\Gamma'_1$ into the point of $I$ lying over the
point $(u_1,
\phi \gamma_1u_1) = (u_1,\theta(\gamma_1)\phi u_1)$
corresponding to the isomorphism $\gamma_1\psi$ of the
pullback of $\eta_1$ to $u_1$ with the pullback of $\eta_2$
to $\phi
\gamma_1u_1$. The morphism $\sigma$ is an isomorphism of
$U_1^\sh\times\Gamma'_1$ with $R^\sh$ in the inverse image of
$C\gen$; it also follows from the fact that the action of\/
$\Gamma'$ on $(\eta_1,U_1)$ is essential that $\sigma$ is
injective. Since the inverse image of $C\gen$ is
scheme-theoretically dense in $R^\sh$ and
$U_1^\sh\times\Gamma_1$ is unramified over $U_1$ we see that
$\sigma$ is an isomorphism of $U_1^\sh\times\Gamma'_1$ with $R^\sh$. It follows
that $R^\sh$ is \'etale 
over $U_1^\sh$; analogously it is \'etale over
$U_2^\sh$. So $R$ is \'etale over $U_1$ and $U_2$ at the
points $u_1$ and
$u_2$; since this holds for all $u_1$ and $u_2$ mapping to
the same point of $C$ the conclusion follows.\qed

Compatibility of charts is stable under base change:

\begin{proposition} \label{Prp:compatibility-pullback} 
\begin{enumerate} 
\item Let $(U_1,\eta_1,\Gamma_1)$ and 
$(U_2,\eta_2,\Gamma_2)$ be two compatible charts for a
generic object $\xi$ on
$C\to S$. If $S'\to S$ is an arbitrary morphism, then
$$(S'\times_S
U_1,\eta_1',\Gamma_1)$$ and
$$( S'\times_S U_2,\eta_2',\Gamma_2),$$
where $\eta'_1$ and $\eta'_2$ are the pullbacks of $\eta_1$
and $\eta_1$ to $S'\times_S U_1$ and $S'\times_S U_2$,
are compatible charts for the pullback of $\xi$ to
$(S'\times_S C\to S')\gen$.
\item If $S' \to S$ is \'etale and surjective, then the converse holds.
\end{enumerate} 
\end{proposition}

The proof is immediate from Proposition \ref{Prp:compatibility-local}. 

\subsection{The product chart}
Given two compatible charts $(U_1,\eta_1, \Gamma_1)$ and
$(U_2,\eta_2, \Gamma_2)$, set $\eta = {\rm pr}_1^* \eta_1$
in $\cM(R)$. There is an
action of\/ $\Gamma$, coming from pulling back the action of\/
$\Gamma_1$ on $\eta_1$; also the tautological isomorphism
$\alpha\colon  {\rm pr}_1^* \eta_1 \simeq \eta_2$ induces an
action of\/ $\Gamma_2$ on $\eta$. These two actions commute,
and therefore define an action of\/
$\Gamma_1\times\Gamma_2$ on $\eta$, compatible with the
action of\/ $\Gamma_1\times\Gamma_2$ on $\rho$. Also, $R$
has a structure of an $n$-marked curve, by defining
$\Sigma^R_i$ to be the inverse image of $\Sigma^{U_1}_i$, and
the map $R \to C$ is strict. Then
$$(R,\eta,\Gamma_1\times\Gamma_2)$$
is a chart, called {\it the product chart}. It is compatible
with both of the original charts.

\subsection{Atlases and twisted objects}

\begin{definition} Fix two non-negative integers $g$ and
$n$. An\/ {\em $n$-pointed twisted object  $(\xi, C \to S,\cA)$ of genus $g$}
consists of 
\begin{enumerate}
\item a  proper, $n$-pointed curve  $C\to S$  of finite presentation,
with geometrically connected fibers of genus $g$;
\item  a generic object $\xi$ over 
$C\to S$; and 
\item  a collection
$\cA = \{(U_\alpha, \eta_\alpha,\Gamma_\alpha)\}$ of
mutually compatible charts, such that the images of the
$U_\alpha$ cover $C$.
\end{enumerate}

A collection of charts $\cA$ as in (3) is called an\/ {\em atlas}.

A twisted object is called\/ {\em balanced} if
each chart in its atlas is balanced (Definition \ref{Def:balanced-chart}).
\end{definition}

\begin{lemma} If two charts for a twisted object are
compatible with all the charts in an atlas, they are mutually
compatible.

Furthermore, if the twisted object is balanced, then any chart which
is compatible with every chart of the atlas is balanced.
\end{lemma} 
\proof Both statements are immediate from the local characterisation of
compatibility in Proposition  \ref{Prp:compatibility-local}.

\begin{remark} The lemma above allows one to define a
twisted object using a maximal atlas, if one prefers.
\end{remark}

\begin{definition}  A {\em morphism of
twisted objects} $(\xi, C \to S,\cA)$ to $(\xi', C' \to
S',\cA')$  consists of a cartesian diagram
$$
\begin{array}{ccc} C& \stackrel{ f}{\longrightarrow} & C'\\
\dar&&\dar\\
S& \longrightarrow & S'\end{array}$$
and an arrow $\xi \to \xi'$ lying over the restriction $f
\rest{C\gen} \colon  C\gen \to C'\gen$, with the property that
the pullback of the charts in $\cA'$
are all compatible with the charts in $\cA$.

The composition of morphisms of twisted objects is defined to be the one
induced by composition of morphisms of 
curves. 
\end{definition}

Let $(\xi, C \to S,\cA)$ be a twisted object, and $T \to S$ a
morphism. Then, using Proposition \ref{Prp:compatibility-pullback} one can
define 
the {\em pullback} of $(\xi, C \to S,\cA)$ to $T$ in the obvious way. 

\subsection{Stability}

\begin{lemma} Let  $(\xi, C \to S,\cA)$ be a
twisted object. Then the morphism $C\gen\to \bM$ induced
by  $\xi$ extends uniquely to
a morphism $C\to \bM$.
\end{lemma}

\proof The unicity is clear from the fact that $\bM$ is separated and
$C\gen$ 
is scheme-theoretically dense in $C$. To prove the existence of an
extension is a local question in the \'etale topology; but if $\cA = \{(
U_\alpha,\eta_\alpha,\Gamma_\alpha)\}$ then  the objects
$\eta_\alpha$ 
induce morphisms $U_\alpha\to \bM$, which are $\Gamma_\alpha$-equivariant,
yielding 
morphisms $U_\alpha/\Gamma_\alpha\to \bM$. These morphisms are
extensions of the pullback
to 
$(U_\alpha)\gen/\Gamma_\alpha$ of the morphism $C\gen\to \bM$. Therefore
they descend to 
$C$.\endproof 

We can now define the main object of this section:

\begin{definition}\label{Def:stable-twisted-object} A twisted object is\/ {\em
stable} if 
the associated map $C \to \bM$ is Kontsevich stable.
\end{definition}

\subsection{The stack of stable twisted objects}

Fix an ample line bundle ${\cO_{\bM}(1)}$ over $\bM$. We
define a category $\TSM = \KO{g}{n}{\cM/\bbS}{d}$ as follows. The objects
are stable twisted objects $(\xi, C \to S,\cA)$, where $C \to S$
is a nodal $n$-pointed curve of genus $g$, such that for the
associated morphism $f \colon  C \to \bM$ the degree of
the line bundle $f^* {\cO_{\bM}}(1)$ is $d$. The arrows are
morphism of twisted objects.

As stated in Theorem \ref{Th:stable-maps}, this category is a proper algebraic 
stack  which is relatively of Deligne--Mumford
type over $\SM$, admitting a projective coarse
moduli sapce 
$\tsm$. The proof of the theorem will begin in section
\ref{Sec:proof}. 

 We shall also consider the full
subcategory $\TSMB$ of {\em balanced} twisted objects.
It will be shown in Proposition \ref{Prp:balanced-is-clopen} that this is an
open and closed substack in  $\TSM$, whose moduli space is open and
closed in $\tsm$.

\section{Twisted curves and twisted stable maps}

In this section we give a stack-theoretic description of twisted
objects. 

\subsection{Nodal stacks}
Let $S$ be a scheme over $\bbS$. Consider a proper, flat, tame
Deligne--Mumford stack
${\cC}\to S$ of finite presentation, such that its fibers
are purely one-dimensional and geometrically connected, with at most 
nodal singularities. Call $C$ the moduli space of ${\cC}$;
By \cite{Keel-Mori} this exists as an algebraic space.

\begin{proposition} The morphism $C \to S$ is a proper flat
 nodal curve of finite presentation, with geometrically
connected fibers.
\end{proposition}

\proof First of all let us show that $C$ is flat over $S$. We
may assume that $S$ is affine; let $R$ be its coordinate ring. Fix a 
geometric point $c_0 \to C$, and call
$C^\sh$ the strict henselization of $C$ at $c_0$. Let $U$ be an
\'etale cover of ${\cC}$, and $u_0$ a geometric point of $U$
lying over $c_0$; denote by $U^\sh$ the strict henselization of
$U$ at $u_0$. If $\Gamma$ is the automorphism group of the
object of ${\cC}$ corresponding to $u_0$, then $\Gamma$
acts on $U^\sh$, and $C^\sh$ is the quotient $U^\sh/ \Gamma$.
Since $\cC$ is tame, the order of $\Gamma$ is prime to the residue 
charactersitic of $u_0$, therefore the coordinate ring of
$C^\sh$ is a direct summand, as an
$R$-module,  of the coordiante ring of $U^\sh$, so it is flat over $R$. 

The fact that the fibers are nodal follows from the fact
that, over an algebraically closed field, the quotient of a
nodal curve by a group action is again a nodal
curve. Properness is clear; the fact that the morphism
${\cC} \to C$ is surjective implies that the
fibers are geometrically connected.\qed

\begin{definition}
A {\em twisted nodal $n$-pointed curve over $S$} is a  diagram
$$\begin{array}{ccc} \Sigma_i^{\cC} & \subset & \cC  \\  &\searrow & \dar  \\ &
 & C   
  \\ &  & \dar  \\ &&S  
\end{array}$$
where
\begin{enumerate}
\item $\cC$ is a tame Deligne-Mumford stack, proper over $S$, and
\'etale locally is a nodal curve 
over $S$;
\item $\Sigma_i^{\cC}  \subset  \cC$ are disjoint closed substacks in the
smooth locus 
of $\cC \to S$;  
\item $\Sigma_i^{\cC} \to S$ are \'etale gerbes; 
\item the morphism $\cC \to C$ exhibit $C$ as the coarse moduli scheme of 
$\cC$; and
\item $\cC \to C$ is an isomorphism over $C\gen$.
\end{enumerate}
\end{definition}

Note that if we let $\Sigma_i^C$ be the coarse moduli spaces of
$\Sigma_i^{\cC}$, then, since $\cC$ is tame, the schemes  
$\Sigma_i^C$ embed in $C$ (they are the images of $\Sigma_i^{\cC}$), and $C$
becomes a 
usual nodal pointed curve. We say that $\cC \to S$ is a twisted pointed curve
{\em of genus $g$,} if $C \to 
S$ is a pointed curve of genus $g$.

\subsection{Morphisms of twisted $n$-pointed nodal curves}

\begin{definition} Let $\cC \to S$ and $\cC' \to S'$ be twisted $n$-pointed
nodal 
curves. A {\em $1$-morphism} (or just a {\em morphism})  $F \colon  \cC \to
\cC'$ is 
a  {\em cartesian} diagram 
$$\begin{array}{ccc}
\cC  & \stackrel{F}{\to} & \cC' \\
\dar &                   & \dar \\
S    & \stackrel{f}{\to} & S' 
\end{array}$$
such that $F^{-1}\Sigma_i^{\cC'} = \Sigma_i^{\cC}$.

If $F,F_1 \colon  \cC \to \cC'$ are morphisms, then we define a 2-morphism $F
\to 
F_1$ to be a base preserving natural transformation. (This implies
that it is an isomorphism.)
\end{definition}

In this way, twisted pointed curves form
a 2-category. However, we have the following:

\begin{proposition}[\cite{A-V:stable-maps}]\label{Prp:category} The 2-category
of twisted pointed 
curves is equivalent (in the lax sense, \cite{K-S}) to a category. 
\end{proposition}

We call the resulting category {\em the category of twisted pointed curves}.
A morphism in this category is an isomorphism class of 1-morphisms in the
2-category of twisted pointed curves.

\proof Since all 2-morphisms are invertible, this claim is the
same as saying that a 1-arrow in the category cannot have
nontrivial automorphisms. The point here is that the stack
${\cC}$ has a dense open {\em representable} substack, which is sufficient by
the following lemma.

\begin{lemma} Let $F \colon   {\cX} \to {\cY}$ be a representable
morphism of Deligne--Mumford stacks over a scheme $S$. Assume
that there exists a dense open representable substack (i.e. an algebraic space)
$U 
\subseteq {\cX}$ and an open representable substack $V
\subseteq {\cY}$  such that $F$ maps $U$ into $V$. Further assume that the
diagonal 
${\cY} \to {\cY}\times_S {\cY}$ is separated. Then
any automorphism of $F$ is trivial.
\end{lemma}

\proof 
\begin{enumerate}
\item
First note that the lemma holds if $ {\cX} = X$ is an algebraic space: 
denote by $\eta$ the 
object of ${\cY}$ over $X$ corresponding to $F$.
The fact that  the diagonal ${\cY} \to {\cY}\times_S {\cY}$ is separated
implies that the isomorphism scheme $\Isom_X(\eta,\eta)\to X$ is
separated. Since $V$ is representable, we have that $\Isom_U(\eta|_U,\eta|_U)\to
U$ is an isomorphism. Thus the unique section over the given open set $U$ has
at most one extension to $X$, which gives the assertion in this case. 

We
will now  reduce the general case to this one by descent. We start with some
observations. 

\item
Let $\alpha$ be an
automorphism of $F$; for each object $\xi$ of ${\cX}$ over
a scheme $X$ we are given an automorphism $\alpha_\xi$ of
$F(\xi)$ over $X$, satisfying
the usual condition for being a
natural transformation. We are going to need the following
two facts. 

\begin{enumerate}
\item If $\xi \to \zeta$ is a morphism in
${\cX}$, then if $\alpha_\zeta$ is
trivial then also $\alpha_\xi$ is trivial. 
This follows immediately from the fact that ${\cX}$ is
a category fibered in groupoids (see
\ref{Sec:stack-criteria}-(\ref{It:pullbacks-unique})).
\item If $X' \to X$ is an
\'etale surjective map of schemes and $\xi'$ is the pullback
of $\xi$ to $X'$, then  if $\alpha_{\xi'}$ is trivial then also
$\alpha_\xi$ is trivial. This follows from the fact that the isomorphism
functors of $\cY$ are sheaves in the \'etale topology (see
\ref{Sec:stack-criteria}-(\ref{It:Isom-sheaf})). 
\end{enumerate}

\item
 Let $X \to {\cX}$ be an \'etale
cover, and call $\xi$ the corresponding object of ${\cX}$
over $X$. The restriction of $\alpha_\xi$ to the open
subscheme $\tilde{U} = U \times_{\cX} X \subseteq X$ is trivial,
and $\tilde{U}$ is scheme-theoretically dense in $X$; applying this lemma in
the case of algebraic spaces it follows that
$\alpha_\xi$ is trivial.

\item
Take an arbitrary object $\tau$ of ${\cX}$ over a scheme
$T$; then there is an \'etale cover $T' \to T$ such that the
pullback $\tau'$ of $\tau$ to $T'$ admits a morphism $\tau
\to \xi$. Applying the two facts mentioned above, we get
the result.\end{enumerate}\endproof

\subsection{Twisted stable maps into a stack} As before, we consider a
proper tame Deligne-Mumford stack $\cM$ admitting a projective coarse
moduli scheme
$\bM$. We fix an ample invertible sheaf on $\bM$.

\begin{definition}\label{Def:twisted-stable-map} A {\em twisted stable
$n$-pointed map of genus $g$ and degree $d$ over $S$}
$$(\cC \to S, \Sigma_i^{\cC}\subset \cC, f\colon  \cC \to \cM)$$   consists
of a commutative diagram
$$\begin{array}{ccc} \cC &\to& \cM \\
		\dar & & \dar \\
		C & \to & \bM \\
\dar &&\\ S &&
\end{array}$$
 along with $n$ closed substacks $\Sigma_i^{\cC}\subset \cC$,
satisfying:
\begin{enumerate}
\item $\cC \to C \to S$ along with $\Sigma_i^{\cC}$  is a twisted nodal
$n$-pointed curve  over $S$; 
\item the morphism $\cC \to \cM$ is representable; and
\item $(C\to S, \Sigma_i^C, f\colon C \to \bM)$ is a stable $n$-pointed map
of degree $d$.   
\end{enumerate}
\end{definition}

\begin{definition} A {\em $1$-morphism} (or just a {\em morphism}) of
twisted stable maps 
$$ G\colon  (\cC \to S, \Sigma_i^{\cC}, f\colon  \cC \to \cM) \to (\cC' \to S',
  \Sigma_i^{\cC'}, f'\colon  \cC' \to \cM)$$
 consists of data $G=(F,\alpha)$, where $F\colon \cC \to \cC'$ is a morphism
of twisted pointed curves, and  $\alpha\colon  f \to f'\circ F$ is an
isomorphism. 
\end{definition}

Twisted stable maps naturally form a 2-category. But by  Proposition
\ref{Prp:category},  this 2-category is equivalent to a category. We
call this category {\em the category of twisted stable maps}. 

\subsection{Equivalence of stable twisted objects and twisted stable
maps}

\begin{theorem} The category of twisted stable maps
is equivalent to the category of stable
twisted objects, via an equivalence which preserves base schemes.
\end{theorem}

\proof {\bf First part:} construction of a functor from stable twisted
objects to the category of twisted stable maps. 

{\sc Step 1: construction of $\cC$ given $(\eta,C,\cA)$.} Consider a
stable twisted object $(\eta, C \to S,\cA)$ with  $\cA=\{(U_\alpha,
\eta_\alpha,\Gamma_\alpha)\}$. For each pair of indices
$(\alpha,\beta)$ let
$R_{\alpha\beta}$ be the transition scheme from $(U_\alpha,
\eta_\alpha,\Gamma_\alpha)$ to $(U_\beta,
\eta_\beta,\Gamma_\beta)$. Let
$U$ be the disjoint union of the $U_\alpha$, and let $R$ be the
disjoint union of  the $R_{\alpha\beta}$. The definition of
$R_{\alpha\beta}$ via isomorphisms (see \ref{Sec:transition scheme}) 
implies that these have the following structure:
\begin{itemize}
\item there are two projection $R\to U$, which are
\'etale;
\item there is a natural diagonal morphism $U\to R$ which sends each
$U_\alpha$ to $R_{\alpha\alpha}$; and
\item there is a product $R\times_U R\to R$, sending each
$R_{\alpha\beta} \times_{U_\beta} R_{\beta\gamma}$ to
$R_{\alpha\gamma}$ via composition of isomorphisms.
\end{itemize}
 These maps give $R\double U$ the structure of an \'etale groupoid,
which defines a quotient Deligne-Mumford stack, which we denote by
$\cC$. This is obviously a  nodal stack on $S$, and its moduli space
is $C$. It is clear that $\cC \to C$ is an isomorphism over $C\gen$. Also
note that, \'etale locally over $C$, the stack $\cC$ is isomorphic to
$[U_\alpha/\Gamma_\alpha]$. 

This construction  depends on the atlas chosen, however we have: 
\begin{lemma}
Let $\cA'$ and $\cA''$ be two compatible atlases on a generic object $(\xi,
C)$.  
Then the stacks $\cC', \cC''$ associated to the corresponding twisted objects
are canonically 
isomorphic.
\end{lemma}
\proof 
In fact, let $\cA$ be the union of
$\cA'$ and $\cA''$. Let $R'\double U'$,
$R''\double U''$ and $R\double U$ the groupoids constructed from these
three atlases, $\cC'$, $\cC''$ and $\cC$ the quotient stacks. There
are obvious embeddings of $R'\double U'$ and $R''\double U''$ into
$R\double U$ inducing isomorphisms of $\cC'$ and $\cC''$ with $\cC$; by
composing the isomorphism $\cC'\simeq \cC$ with the inverse of
$\cC''\simeq
\cC$ we obtain the desired canonical isomorphism $\cC'\simeq \cC''$.\qed

{\sc Step 2: construction of $\Sigma_i^{\cC}$.} Since, given two
indices $\alpha, \beta$, the inverse image of
$\Sigma^{U_\alpha}_i$ in $R_{\alpha\beta}$ coincides with the inverse
image of $\Sigma^{U_\beta}_i$, the collection of the
$\Sigma^{U_\alpha}_i$ defines a closed substack $\Sigma^{\cC}_i$ of
${\cC}$. Since the $\Sigma^{U_\alpha}_i$ are \'etale over
$S$, it follows that $\Sigma^{\cC}_i$ is \'etale over $S$;
furthermore, the moduli space of $\Sigma^{\cC}_i$ is $\Sigma^C_i$, so
for any algebraically closed field
$\Omega$ the induced functor $\Sigma^{\cC}_i(\Spec \Omega)
\to \Sigma^C_i(\Spec \Omega)$ induces a bijection of the set of
isomorphism classes in $\Sigma^{\cC}_i(\Spec
\Omega)$ and $\Sigma^C_i(\Spec \Omega)$. This means that
$\Sigma^{\cC}_i$ is an \'etale gerbe over $S$, and ${\cC} \to S$ has
the structure of a twisted $n$-pointed curve. 

If we start from a different, but compatible, atlas, the isomorphism
between the two twisted curves constructed above preserves these
markings.

{\sc Step 3: construction of $Ob({\cC}) \to Ob(\cM)$.} Putting
together the  objects $\eta_\alpha$ we have an object
$\eta$ on $U$, and the tautological isomorphism between the two
pullbacks of $\eta_\alpha$ and $\eta_\beta$ to
$R_{\alpha\beta}$ yield an isomorphism $\phi$ of the two pullbacks of
$\eta$ to $R$. Fix a morphism $T \to {\cC}$, where $T$ is a scheme,
and  set $T' = T 
\times_{\cC}U$, $T'' = T
\times_ {\cC}R$. Denote by $\zeta'$ the pullback of
$\eta$ to $T'$, then $\phi$ induces an isomorphism of the two
pullbacks of $\zeta'$ to $T''$ satisfying the cocycle condition. By
descent  for objects of  $\cM$ (see
\ref{Sec:stack-criteria}-(\ref{It:etale-descent})), this isomorphism   
allows  $\zeta'$ to descend to  an object $\zeta$ of ${\cM}(T)$, or
equivalently a morphism
$T \to \cM$. Thus we have associated to an object of $\cC$ an object
of $\cM$.

{\sc Step 4: construction of $Mor({\cC}) \to Mor(\cM)$.} We need to
associate, for each arrow in $\cC$, an arrow in $\cM$ in such a manner
that we get a morphism of stacks  $\cC \to \cM$. 

Consider an arrow in $\cC$ over a morphism of schemes $T_1 \to T_2$,
which corresponds to a commutative diagram
$$\begin{array}{ccccc} T_1 & &\lrar& & T_2 \\
 &\searrow && \swarrow &\\
 &  & \cC  &&
\end{array}$$ Set
$T'_i = T _i\times_{\cC}U$, $T''_i = T_i\times_ {\cC}R$. Denote by
$\zeta'_i$ the pullback of $\eta$ to $T'_i$, $\zeta_i$ the object
obtained by descent to
$T_i$; there is an arrow $\zeta'_1 \to \zeta'_2$, together with
descent data, inducing an arrow $\zeta_1 \to\zeta_2$.

Thus we have morphism of stacks ${\cC} \to \cM$.

If we start from two different but compatible atlases $\cA'$ and
$\cA''$ on the same generic object, we obtain two morphisms $\cC' \to \cM$ and
$\cC'' \to \cM$. It 
is easily seen that there is a canonical isomorphism between this
morphism $\cC'
\to \cM$ and the morphism $\cC' \to \cM$ obtained by composing the
isomorphism
$\cC' \simeq \cC''$ constructed above with the morphism $\cC'' \to
\cM$.

{\sc Step 5: the morphism ${\cC} \to \cM$ is representable.}
 This is a consequence of the fact that the action of the finite
groups appearing in the charts of the atlas is essential, because of
the following Lemma, which is probably well-known.

\begin{lemma} Let $g\colon \cG\to\cF$ be a morphism of Deligne Mumford
stacks. The following two conditions are equivalent:
\begin{enumerate}
\item The morphism  $g\colon \cG\to\cF$ is representable.
\item For any algebraically closed field $k$ and any $\xi\in \cG(k)$,
the natural group homomorphism $\Aut (\xi) \to \Aut(g(\xi))$ is a
monomorphism.
\end{enumerate}
\end{lemma} 
\proof By definition, the morphism $g\colon \cG\to\cF$  is representable if
and only if  the following condition holds: 
\begin{itemize}
\item For any algebraic space $X$ and any morphism $f\colon X \to \cF$, the
stack
$\cY:=\cG \times_\cF X$ is equivalent to an algebraic space.
\end{itemize} For fixed $f\colon X \to \cF$ latter condition is equivalent
to the following:
\begin{itemize}
\item The diagonal morphism $\Delta\colon \cY \to \cY\times \cY$ is a
monomorphism.
\end{itemize} This means:
\begin{itemize}
\item Given an algebraically closed field $k$, and an element $q\in 
\cY\times
\cY(k)$ and two elements $p_i\in \cY(k)$, with isomorphisms
$\beta_i\colon \Delta(p_i) 
\to q$, there exists a unique isomorphism $\phi\colon p_1 \to p_2$ such that
$\beta_1= \beta_2\circ\Delta(\phi)$. 
\end{itemize} We can write $$q = \left((\g_1,x_1, \alpha_1),(\g_2,x_2,
\alpha_2)\right)$$ where $\g_i\in \cG(k), x_i\in X(k)$, and $
\alpha_i\colon g(\g_i)
\stackrel{\sim}{\to}  f(x_i)$. Similarly $p_i=(\g'_i,x'_i, \alpha'_i)$
as above. The existence of
$\beta_i$ implies $x_1 = x_2 = x'_1 = x'_2$. Also, composing with  the
given isomorphisms  we may reduce to the case that, in fact,  the same
is true for the $\g_i,
\g_i', \alpha_i$ and $\alpha_i'$. Thus the condition above is
equivalent to saying that for any 
$\g\in G(k)$ there is a unique $\beta\in \Aut(\g)$ such that $g(\beta)
= id_{g(\g)}$, which is what we wanted. \qed

Now, let $q \in \cM(\Spec k)$, where $k$ is an algebraically closed
field, and lift $q$ to a geometric point $p \in U(\Spec k)$; let
$\alpha$ such that $p \in U_\alpha$. The automorphism group of $q$ in
$\cC$ is the stabilizer of $p$ inside $\Gamma_\alpha$; the fact that
this acts faithfully on the fiber of
$\eta_\alpha$ over $p$ means exactly that this automorphism group
injects inside the automorphism group of the image of $q$ in $\cM$.

{\sc Step 6: stability.}

The map $C \to S$ induced by this morphism ${\cC} \to S$ is the map
induced by the twisted object $(\xi, C \to S,\cA)$, and so is stable; this
way we get a twisted stable map ${\cC} \to {\cM}$.

{\sc Step 7: construction of the morphism of twisted stable maps
induced by a morphism of twisted objects.} 

Let $f \colon  (\xi_1, C_1 \to S_1,\cA_1) \to (\xi_2, C_2 \to S_2,\cA_2)$ be a
morphism of twisted objects, and let $\zeta_1\colon  \cC_1 \to \cM$ and
$\zeta_2\colon  \cC_2
\to \cM$ the corresponding twisted stable maps. Denote also by  
$\zeta_1'\colon  
\cC_1' \to \cM$ the twisted stable map associated to the twisted object
$(\xi_1, C_1 \to S_1,f^* \cA_2)$. Since, as we have 
seen, there is a canonical isomorphism between $\zeta_1$ and $\zeta_1'$, we are
reduced to the case $\cA_1 = f^* \cA_2$. In this case there is an morphism of
groupoids 
$R_1\double U_1$ to $R_1\double U_1$, which induces a morphism of
stacks $\cC_1 \to \cC_2$.

This completes the definition of the functor from twisted objects to
twisted stable maps.

{\bf Second part:} construction of a functor 
 from twisted stable maps to twisted objects.

{\sc Step 1: construction of $(\eta,C,\cA)$.} Let us take a twisted stable
map ${\cC}
\to {\cM}$, with moduli space $C \to S$. Since $C\gen$ is isomorphic
to it inverse image in ${\cC}$, the map
${\cC} \to {\cM}$ induces a map $C\gen \to {\cM}$, and correspondingly
a generic object $\xi$ on $C\gen$.

By lemma \ref{Lem:locally-quotient}, there is an \'etale covering
$\{C_\alpha \to C\}$, such that for each
$\alpha$ there is a scheme $U_\sigma$ and a finite group
$\Gamma_\alpha$ acting on $U_\alpha$, with the property that the
pullback
${\cC}\times_CC _\alpha$ is isomorphic to the stack-theoretic quotient
$[U_\alpha/\Gamma_\alpha]$.

For each $i= 1, \ldots, n$ we define $\Sigma^{U_\alpha}_i$ as the
pullback of $\Sigma^{\cC}_i$ to $U_\alpha$; then the action of the
groups $\Gamma_\alpha$ preserves the subschemes
$\Sigma^{U_\alpha}_i \subseteq C$. These subschemes are obviously
disjoint and contained in the smooth locus of the projection
$U_\alpha \to S$, so $U_\alpha \to S$ is an $n$-marked curve.

If we call $\eta_\alpha$ the object of ${\cM}(U_\alpha)$ corresponding
to the morphism $U_\alpha \to {\cM}$, the action of $\Gamma_\alpha$ on
$U_\alpha$ lifts to an action of
$\Gamma_\alpha$ on $\eta_\alpha$. We get charts
$(U_\alpha, \eta_\alpha, \Gamma_\alpha)$ for the generic object
$\xi$. The fiber product
$U_\alpha \times_ {\cC} U_\beta$ is \'etale over
$U_\alpha$ and $U_\beta$. Since the given morphism ${\cC}
\to {\cM}$ is representable and separated we have that
$U_\alpha \times_ {\cC} U_\beta$ is a closed subspace of
$$ U_\alpha \times_ {\cM} U_\beta = {\rm Isom}_{U_\alpha
\times_S U_\beta}({\rm pr}_\alpha^*\eta_\alpha, {\rm
pr}_\beta^*\eta_\beta),
$$ where ${\rm pr}_\alpha$ and ${\rm pr}_\beta$ denote the projections
of $U_\alpha \times_S U_\beta$ onto the two factors. This implies that
$U_\alpha \times_ {\cC} U_\beta$ is the transition scheme of the two
charts, which are therefore compatible. It follows that we have
defined an atlas on $C \to S$, and so a twisted object. The induced
morphism
$C \to {\bM}$ is the one induced by ${\cC} \to {\cM}$, so it is
stable. 

If we start from a diffent covering of $\cC$ by
$[U_\alpha/\Gamma_\alpha]$ we get a canonically isomorphic twisted
object. In fact, we essentially get the same twisted object with a
different atlas.

{\sc Step 2: construction of morphisms.} Given a morphism from $\cC_1
\to \cM$ to $\cC_2 \to \cM$, we automatically get a morphism $C_1 \to
C_2$. This association is obviously compatible with composition of
morphisms.

Since $\cC_1$ and $\cC_2$ are tame, the formation of their moduli
spaces $C_1$ and $C_2$ commute with base change (Lemma
\ref{Lem:tame-cms-pullback}), so that 
$C_1 = S_1
\times_{S_2}C_2$; this implies that  $\cC_1 $ is isomorphic to $C_1
\times_{C_2} \cC_2$, and the morphism $\cC_1 \to \cM$ is isomorphic to
the composite  $\cC_1 \to \cC_2 \to
\cM$. From this it follows immediately that the pullback of a chart
$(U,\eta,\Gamma)\in \cA_2$, is compatible with all charts in
$\cA_1$.

{\bf Third part:} proving that the two functors  defined above are 
inverse to each other.

First of all start from a twisted map $\cC \to \cM$, and call $(\eta,
C,\cA)$ the associated twisted object. This is obtained by taking a
covering of $\cC$ via $U_\alpha$ as in part~2; if we call
$U$ the disjoint union of the $U_\alpha$, and let $R = U \times_\cC
U$, the twisted curve associated with $(\eta, C,\cA)$ is the quotient of
the groupoid $R\double U$, which is (canonically isomorphic to) $\cC$.
It is straightforward to check that the morphism $\cC \to \cM$
obtained from $(\eta, C,\cA)$ is canonically isomorphic to the given one.

Now start from a twisted object $(\eta, C,\cA)$, and consider
the associated twisted map $\cC \to \cM$. The atlas
$\cA = \{(\eta_\alpha, U_\alpha, \Gamma_\alpha)\}$ yields a covering of
$\cC$ with stacks $[U_\alpha/\Gamma_\alpha]$; the associated twisted
object is canonically isomorphic to $(\eta, C,\cA)$.

\endproof 
\begin{remark}
We draw attention to the following counter-intuitive phenomenon: given a
twisted stable map  $f\colon  \cC \to \cM$, we can look at the automorphisms
in $\Aut_C\cC$. Since $\cC\gen = C\gen$, such an automorphism is a {\em local}
object, dictated by the 
structure of $\cC$ along the special locus of $C$. However, if such an
automorphism comes from an automorphism of  $f\colon  \cC \to \cM$, then it is
determined by the associated automorphism of the {\em generic} object $\xi$ on
$C$.
\end{remark}

\section{The category $\TSM$ is an algebraic stack}\label{Sec:proof}

\subsection{The stack axioms}

\begin{proposition}\label{Prp:tsm-stack} The category $\TSM$ is a
limit-preserving stack, fibered by 
 groupoids over ${\cS}ch/\bbS$. 
\end{proposition}
\proof 
{\sc Condition \ref{Sec:stack-criteria}-(\ref{It:fibered-by-groupoids}):} By
definition  $\TSM$ is  fibered by groupoids  
over ${\cS}ch/\bbS$, since there are pullbacks, and since all the
morphisms of objects are given by fiber diagrams. 

{\sc Condition \rm \ref{Sec:stack-criteria}-(\ref{It:limit-preserving}):} It is
also 
not difficult to see that $\TSM$ is limit 
preserving: given a twisted object $(\xi, C,\cA)$ over $\Spec R$, where $R =
\indlim R_i$, the schemes $C, U_\alpha$ are of
finite 
presentation, and therefore they come from some $\Spec R_i$; since the stack
$\cM$ is limit preserving,  the objects $\eta_\alpha$ come from 
$\Spec R_{i'}$ for some $i'\geq i$. 

{\sc  Condition \rm \ref{Sec:stack-criteria}-(\ref{It:etale-descent}):} We need
to 
show that 
$\TSM$ has effective \'etale 
descent for objects. 

Given a scheme $S$, an
\'etale cover $\{S_\alpha\to S\}$, and twisted objects
$(\xi_\alpha,C_\alpha \to S_\alpha,\cA_\alpha)$
 together with isomorphisms
between the pullbacks of $(\xi_\alpha,C_\alpha \to S_\alpha,\cA_\alpha)$ and
$(\xi_\beta,C_\beta \to S_\beta,\cA_\beta)$ to $S_\alpha\times_S S_\beta$
satisfying the 
cocycle condition, we claim that these descend to a twisted object
$(\xi,C\to S,\cA)$ on $S$. The existence of the projective curve $C\to S$  
 is immediate from the sheaf axiom for $\SM$. Now, any chart
$(U_\alpha,\eta_\alpha,\Gamma_\alpha)$ for $(\xi_\alpha,C_\alpha \to
S_\alpha)$ is also automatically a chart for $(\xi,C\to S)$.
The charts coming from one $(\xi_\alpha,C_\alpha \to S_\alpha)$ are
obviousely compatible. Compatibility 
of the charts 
coming from $(\xi_\alpha,C_\alpha \to S_\alpha)$ and
$(\xi_\beta,C_\beta \to S_\beta)$ follows since they are compatible 
when pulled back to $S_\alpha\times_S S_\beta$.

{\sc Condition \rm \ref{Sec:stack-criteria}-(\ref{It:Isom-sheaf}) and
(\ref{It:Isom-representable}):} The other  sheaf axiom requires each $\Isom$ 
functor to be a sheaf in the 
\'etale topology. Since we would like to show that $\TSM$
is an {\em 
algebraic} stack, 
we might as well  prove that the $\Isom$ functor is
representable:
\begin{proposition}\label{Prp:isom-representable} For any pair of stable
twisted 
objects 
$$\tau_1= (\xi_1, C_1 \to S,\cA_1) \quad \mbox{and}\quad \tau_2= (\xi_2, C_2
\to S,\cA_2)$$ over the same scheme $S$, the 
functor $\Isom_S(\tau_1,\tau_2)$ of isomorphisms of
twisted objects is representable by a separated scheme of finite type over $S$.
\end{proposition}

\proof {\sc Step 1: reduction to the case   $C_1 = C_2$}.
Consider the associated maps $\phi_i\colon  C_i \to {\bM}$;
there is a natural transformation $\Isom_S(\tau_1,\tau_2) \to
\Isom_S(\phi_1,\phi_2)$, where the second functor is the one
associated with the diagonal in the stack $\SM$. Since
the second 
functor is known to be representable by separated schemes of finite
type, it is enough to prove that this natural transformation is
representable. This means that we can assume that
$C_1$ is equal to $C_2$, and reformulate the problem as
follows: if $\tau_1 = (\xi_1, C
\to S,\cA_1)$ and
$\tau_2 = (\xi_2, C
\to S,\cA_2)$ are stable twisted objects over the same nodal curve
$C \to S$, then the functor $I$ of isomorphisms of the two
twisted objects inducing the identity on $C$ is representable
by a separated scheme of finite type over $S$.

 {\sc Step 2: reduction to case  where  $\xi_1 = \xi_2$}.
Consider the scheme $\Isom_{C\gen} (\xi_1, \xi_2)$, which is
finite over $C\gen$; it can be  extended,  to a scheme $J$ {\em finite
over $C$},  as follows. 

Let $\cC_1$ and $\cC_2$ be the twisted curves underlying the twisted stable
maps associated to $\xi_1, \xi_2$. Then the scheme $$\Isom_{C\gen} (\xi_1,
\xi_2) = C\gen\mathop{\times}\limits_{\cM\times C}C\gen$$ is open
inside the {\em stack} $$\cC_1 
\mathop{\times}\limits_{\cM\times C}\cC_2.$$ Let $J$ be the coarse moduli space
of the latter 
stack. Then   $J$ is finite over $C$, and $$\Isom_{C\gen} (\xi_1, \xi_2)\subset
J$$ is  open.

By definition, each morphism of twisted objects $\tau_1\to\tau_2$ induces a
section $C \to J$; this defines a morphism from $I$ to the
Weil restriction ${\rm R}_{C/S}J$, which is quasi-projective over $S$ since $C
\to 
S$ is 
projective and $J\to C$ is 
finite (see \cite{Gr-FGA}). 
It is enough to
prove that 
the morphism $I \to {\rm R}_{C/S}J$ is representable, separated and of finite
type, so we may assume 
that $\xi_1$ is equal to $\xi_2$.

{\sc Step 3:  the case   $\xi_1 = \xi_2$.}  We see that the 
following lemma implies the thesis.

\begin{lemma} \label{Lem:isom-closed-global} Let $(\xi, C \to S)$ be a
generic 
object with 
two atlases
${\cA}_1$ and ${\cA}_2$. Then there exists a closed
subscheme
$S'\subseteq S$ such that given a morphism $T\to S$, the
pullbacks of ${\cA}_1$ and
${\cA}_2$ to $T\times_S X$ are compatible if and only if
$T\to S$ factors through $S'$.
\end{lemma}

\proof We will prove a local version of the fact above.

\begin{lemma} \label{Lem:isom-closed-local} Let 
$(\xi, C \to S)$ be a generic object with two charts $(U, \eta, \Gamma)$ and
$(U', \eta', \Gamma')$. Let 
$s_0$ be a geometric point of
$S$; assume that the fiber of $C$ over $s_0$ has a
unique special point
$c_0$, and that there are two unique geometric
points $u_0\in U$ and $u'_0 \in U'$ over $c_0$. Then there exists 
 an
\'etale neighborhood $\tilde S$ of $s_0$, and a closed subscheme
$S'\subseteq \tilde S$ such that given a morphism $T\to \tilde S$, the
pullbacks of the two charts to $T\times_{\tilde S} X$ are compatible
if and only if $T\to \tilde S$ factors through $S'$.
\end{lemma}

{\em Proof of the global Lemma \ref{Lem:isom-closed-global} given
the local Lemma \ref{Lem:isom-closed-local}.}

First fix a geometric point $s$ of $S$, and let $c_1,\ldots,c_k$ be the
special points above it. For each $c_i$ let $(U_i, \eta_i, \Gamma_i)$ and
$(U_i', \eta_i', \Gamma_i')$ be charts from ${\cA}_1$ and ${\cA}_2$. Refining
the charts if necessary, we may assume that there are unique geometric points
$u_i$ and ${u_i}'$ in $U_i$ and $U_i'$ over $c_i$. By Lemma
\ref{Lem:isom-closed-local} there exists an \'etale neighborhood $\tilde S_i$
of $s$ and a closed subscheme $S'_i \subset \tilde S_i$ for compatibility of 
$(U_i, \eta_i, \Gamma_i)$ and
$(U_i', \eta_i', \Gamma_i')$. We can choose a common refinement $\tilde S_s$ of
$\tilde S_i$ such that $C|_{\tilde S_s\setmin \{s\}}$ is of constant
topological type, and the pullbacks 
of $U_i$ and
$U_i'$ cover all special points over $\tilde S_s$. Then the intersection  $S'_s
= \cap S'_i\times_{\tilde S_i} \tilde S_s$ is a closed subscheme of $\tilde
S_s$ such that 
given a morphism $T\to \tilde S_s$, the 
pullbacks of the two atlases to $T\times_{\tilde S_s} X$ are compatible
if and only if $T\to \tilde S_s$ factors through $S'$.

Now choose a finite number of geometric points $s_i$ and \'etale neighborhoods
$\tilde S_{s_i}$ which cover $S$. Applying Proposition
\ref{Prp:compatibility-pullback}, the union of the images of the closed
subschemes $S'_{s_i}$ is the closed subscheme required. \qed


{\em Proof of the local lemma.} By passing
to the fiber product $(U/ \Gamma)
\times_C (U'/ \Gamma')$ we may assume that $C = U/\Gamma =
U'/\Gamma'$. By refining $C$ we may assume that $C$ and $U$
are affine, and that there exists an invariant
effective Cartier divisor $D$ on
$U$ containing the locus where the projection
$U\to S$ is not smooth, but none of the fibers of $U\to S $. Using the local
criterion for flatness, it is easy to see that such a divisor $D$ is flat.


 By passing to an
\'etale neighborhood of $s_0$ we can split $D$ into a number
of connected components, so that the component containing
$v_0$ is finite over $S$; then by deleting the other
connected components we see that we may assume that $D$ is
finite over $S$.

Consider the transition scheme $R$ of the two charts. There
is a free action of the group $\Gamma'$ over $R$; set $E =
R/ \Gamma'$. The projection $E\to S$ is an isomorphism over
the smooth locus of $C$, and is an isomorphism if and only if
the projection $R\to U$ is \'etale. Take a morphism $T\to S$;
then the transition schemes $R_T$ of the pullbacks of the two
charts to $T$ is the scheme-theoretic closure of the
inverse image of the smooth locus of $C$ in $R \times_S
T$, so $R_T$ is \'etale over $U\times_S T$ if and only if the
projection $(R \times_S T)/ \Gamma'  = E\times_S T\to U
\times_ST$ has a section.

Set $S = \Spec \Lambda$, $U = \Spec A$, $E = \Spec B$,
and call $I$ the ideal of $D$ in $A$. The coordinate ring of the 
complement of $D$ inside $U$ is $A' = \cup_{i=0}^\infty
I^{-i}$, and there is a natural homomorphism $B\to A'$;
given a $\Lambda$-algebra $L$ the coordinate ring of the
quotient of the transition scheme of the pullbacks of the two charts to
$\Spec L$ is the image of $B \otimes_\Lambda L$ in $A'\otimes_A L$, so
$(R \times_S T)/ \Gamma'  = E\times_S T\to U
\times_ST$ has a section if and only if the image of $B
\otimes_\Lambda L$ in $A' \otimes_\Lambda L$ is equal to $A
\otimes _\Lambda L\subseteq A' \otimes_\Lambda L$. Take
a set of generators of $B$ as an $A$-algebra, and call
$e_1,\ldots,e_n$ their images inside $A'/A$; the condition
that the image of $B \otimes_\Lambda L$ be equal to $A
\otimes _\Lambda L$ is equivalent to the condition that the
images of the $e_i$ in $(A'/A)\otimes_\Lambda L$ be zero. Fix
an integer $n$ such that the $e_i$ are all contained in
$I^{-n}/A$; then the sequence
$$
0\longrightarrow I^{-n}/A\longrightarrow A'/A\longrightarrow
\cup_{i=n}^\infty I^{-i}/I^{-n}\longrightarrow 0
$$
is exact and $\cup_{i=n}^\infty I^{-i}/I^{-n} =
I^{-n}\otimes_A(A'/A)$ is flat over $\Lambda$, so the
sequence stays exact after tensoring with $L$. The conclusion
is that
$R_T$ is \'etale over $U\times_S T$ if and only if the images
of the $e_i$ in $(I^{-n}/A)\otimes_\Lambda L$ is zero. But
$I^{-n}/A$ is finite and flat as a $\Lambda$-module, so it is
projective, and can be embedded as a direct summand of a free
$\Lambda$-modules $F$. If $J \subseteq \Lambda$ is the ideal
generated by the coefficients of the $e_i$ with respect
to  a basis of $F$, it is clear that the closed subscheme $S'=
\Spec (\Lambda/J)$ has the desired property.\endproof

This completes the proof of Proposition \ref{Prp:isom-representable}, and
also of Proposition \ref{Prp:tsm-stack}.\endproof

\subsection{Base change}
Artin's criteria for an algebraic stack work over a scheme of finite type over
a  field or an excellent Dedekind domain. Since $\cM\to \bbS$ is of finite
presentation, it is obtained by base change from a scheme of finite type over
$\ZZ$. We need a similar statement of $\TSM$, which follows from the following
result:

\begin{proposition}\label{Prp:base-change}
Let $\bbS_1 \to \bbS$ be a morphism of noetherian schemes and let $\cM \to
\bbS$ be as in the main theorem. Then $\KO{g}{n}{(\cM \times_\bbS
\bbS_1)/\bbS_1}{d} =  \TSM \times_\bbS 
\bbS_1$.   
\end{proposition}
\proof Since $\cM$ is tame, by Lemma \ref{Lem:tame-cms-pullback}  we have that
$\bM\times_\bbS 
\bbS_1$ is the moduli space of $\cM\times_\bbS
\bbS_1$. From this the result is immediate. \endproof 


\subsection{Deformations and obstructions}
So far we have seen that $\TSM$ is a limit preserving
stack, with representable diagonal. In order to show that it is
algebraic, we need to  produce an \'etale covering by a scheme $V \to
\TSM$.  We follow Artin's method, in which one starts from a deformation and
obstruction theory, one constructs 
formal deformation spaces, and one shows they are  algebraizable. We start by
constructing a deformation and 
obstruction theory for $\TSM$. Here it is convenient to
realize $\TSM$ as the category of twisted stable maps, rather than stable
twisted objects.

\begin{proposition} The category $\TSM$ has a deformation
and obstruction theory satisfying  conditions
\ref{Sec:stack-criteria}-(\ref{It:def-obs}) (i) through (iii). 
\end{proposition}

Let $A_0$ be a {\em reduced} ${\cO}_ {\bbS}$-algebra of finite type,
$({\cC}_0 \to \Spec A_0, \Sigma_0,{\cC}_0 \to {\cM})$ a twisted stable map over
$A_0$. For 
any 
$A_0$-module
$I$ of finite type, let us call ${\rm D}(I)$ the set of isomorphism classes
of twisted stable maps
$({\cC} \to  \Spec (A_0+I), \Sigma,{\cC} \to  {\cM} )$, with a given
isomorphism 
of the 
restriction of ${\cC} \to {\cM}$ to $A_0$ with ${\cC}_0 \to {\cM}$. Recall that
the set ${\rm D}(I)$ has a natural structure of 
$A_0$-module. Denote by ${\rm L}$ the cotangent complex ${\rm
L}_{{\cC}_0/ {\cM}\times_S \Spec A_0}$ (see \cite{Illusie}, II, 1.2.7,
\cite{L-MB}, 
17.3) and by 
${\cN}$ the normal sheaf of $\Sigma_0$ in ${\cC}_0$.

\begin{lemma} There exists a canonical exact sequence ${\rm E}(I)$ of
$A_0$-modules
$$
\Hom_{{\cC}_0}({\rm L}, I \otimes_{A_0} {\cO}_{{\cC}_0}) \to
\opH^0({\cC}_0, I\otimes {\cN}) \to {\rm D}(I) \to \Ext^1({\rm L}, I
\otimes_{A_0} {\cO}_{{\cC}_0}) \to 0
$$
which is functorial in $I$.
\end{lemma}

\proof By \cite{Illusie}, III, Proposition 2.1.2.3, the group ${\rm D}'(I)
\eqdef 
\Ext^1({\rm L}, I
\otimes_{A_0} {\cO}_{{\cC}_0})$ classifies extensions ${\cC} \to
{\cM}$ of the twisted map ${\cC}_0 \to {\cM}$, with the markings
ignored. The natural map ${\rm D}(I) \to {\rm D}'(I)$ is surjective, because
$\Sigma_0$ is unobstructed inside ${\cC}_0$, as $\opH^1( {\cC}_0,
{\cN}) = 0$. The group $\opH^0( {\cC}_0, I \otimes_{A_0}{\cN})$ classifies
extensions 
of 
$\Sigma_0$ inside the trivial extension ${\cC}_0\times_{\Spec A_0}
\Spec(A_0 + I)$, so $\opH^0( {\cC}_0, I \otimes_{A_0}{\cN})$ surjects onto the
kernel 
of the map ${\rm D}(I) \to {\rm D}'(I)$. Finally, $\Hom_{{\cC}_0}({\rm
L}, I \otimes_{A_0} {\cO}_{{\cC}_0})$ is the group of
${\cO}_{\cM}$-derivations from ${\cO}_{{\cC}_0}$ into $I \otimes_{A_0} 
{\cO}_{{\cC}_0}$, therefore it is the group of infinitesimal automorphisms
of ${\cC}_0\times_{\Spec A_0}
\Spec(A_0 + I)$ over ${\cM}$ fixing ${\cC}_0$, so it surjects onto the
kernel of $\opH^0( {\cC}_0, I \otimes_{A_0}{\cN}) \to {\rm D}(I)$. \endproof 

\begin{lemma} Let $A' \to A \to A_0$ be infinitesimal
extensions as in \cite{Artin}, so that $I = \ker ( A' \to A)$ is an
$A_0$-module, and let  $({\cC} \to \Spec A, \Sigma, {\cC} \to {\cM})$ be a
twisted pointed map. The obstruction to lifting $({\cC} \to \Spec A, \Sigma,
{\cC} \to {\cM})$ to $A'$ lies 
in the group $\Ext^2_{{\cO}_{\cC}}({\rm L}_{{\cC}/ {\cM}}, I
\otimes_A {\cO}_{\cC})$.
\end{lemma}
\proof  If there is an extension ${\cC}' \to {\cM} $ of ${\cC}
\to {\cM}$, then $\Sigma$  lifts to a subgerbe of ${\cC}$, as
indicated above. So we see that the obstruction to extending the twisted
stable map coincides with the obstruction to extending ${\cC}\to \cM$.
According to \cite{Illusie}, III, Proposition 2.1.2.3 this obstruction is
an element 
of $\Ext^2_{{\cO}_{\cC}}({\rm L}_{{\cC}/ {\cM}}, I \otimes_A
{\cO}_{\cC})$. \qed

\begin{lemma} Let $A \to A_0$ be an infinitesimal extension of $A_0$,
and let ${\cC} \to \Spec A$ be a tame proper Deligne--Mumford stack. Let
${\cL}$ be a complex  bounded above,  of sheaves of ${\cO}_{\cC}$ modules, with
coherent  
cohomology. For each \'etale morphism $A \to B$ and each finite $B$-module
$J$, set
$$
{\rm F}_B(J) = \Ext^n_{{\cO}_{{\cC}_B}}({\cL}\otimes_A B, J
\otimes_A {\cO}_{\cC}).
$$
Let $I$ be a finite $A_0$-module. Then

\begin{enumerate}
\item  If $A \to B$ is \'etale, then
$$
{\rm F}_B(I \otimes_{A_0} B_0) = {\rm F}_A(I) \otimes_{A_0}B_0.
$$

\item  If $\m$ is a maximal ideal of $A$, then
$$
{\rm F}_{A_0}(I)\otimes _{A_0} \projlim(A/ \m^k) = \projlim {\rm
F}_{A_0}(I/\m^kI).
$$

\item There is an open dense subset $U \subseteq \Spec A_0$ such that
$$
{\rm F}_{A_0}(I)\otimes _{A_0} k(p) = {\rm F}_{A_0}\left(I\otimes _{A_0}
k(p)\right)
$$
for all $p \in U$.

\end{enumerate}

\end{lemma}

\proof Let ${\cK} = {\rm R}\cH om_{{\cO}_{\cC}}({\cL},
{\cO}_{{\cC}_0})$. Since ${\cL}$ is bounded above, and locally 
quasi-isomorphic to a complex of locally free sheaves, we have
$$
{\rm R}{\cH}om_{{\cO}_{\cC}}({\cL}, {\cO}_{{\cC}_0} \otimes_{A_0}I)
= {\cK}\tototi_{A_0} I. 
$$

By lemma \ref{Lem:cms-exact}, the hypercohomology of the
 complex
${\cK}\tototi_{A_0} I$ is isomorphic to the hypercohomology of
the complex of sheaves $\pi_*({\cK}\tototi_{A_0} I) =
(\pi_*{\cK})\tototi_{A_0} I$ on $C_0\subseteq C$.

Statement 1 follows from the fact that formation of hypercohomology commutes
with flat base change, for separated algebraic spaces; the case of a single
sheaf is in \cite{L-MB}, 13.1.9. From this one can deduce the general
case by the usual spectral sequence argument.

Statement 2 can be deduced from the theorem on formal functions (see
\cite{Knutson}, V, 3.1). Note that  each term in the complexes
on the right is finite over the artinian ring $A_0/m^k$, and therefore the
complexes and their 
cohomologies satisfy the Mittag--Leffler condition. A standard spectral
sequence argument (see \cite{G-EGA3}, 0, 13.2.3) gives the statement.

Statement 3 is proved as follows. Since $A_0$ is reduced, by generic
flatness we can localize, and assume that the cohomology sheaves of
$\pi_*{\cK}\tototi_{A_0} I$ are flat over $A_0$. By localizing further,
we can also assume that the hypercohomology groups of $\pi_*{\cK}\tototi_{A_0}
I$ are projective over $A_0$. Then the statement follows 
from the standard base change theorem (see \cite{SGA6}, IV, 3.10 and proof of
\cite{L-MB}, 13.1.9).
 \endproof

%
%

We now come to our proof of the proposition. 

Condition (i) follows form (1) in the  lemma: for the obstructions it is
immediate, and for the deformations it follows using  the exact sequence ${\rm 
E}(I)$. 

For (ii), first notice that $\projlim {\rm E}(I/\m^k I)$ is still exact,
because each of the terms of ${\rm E}(I/\m^k I)$ is an artinian
$A_0$-module, and so they satisfy the Mittag--Leffler condition. The
sequence ${\rm E}(I) \otimes_{A_0} \widehat A_0$ is also exact, because
$\widehat A_0$ is flat over $A_0$, so the conclusion follows from (2) in the
lemma above, together with the five
lemma.

For (iii), note that there is a dense open set $U \subseteq \Spec A_0$ such
that ${\rm E}(I) \otimes_{A_0} k(p)$ is exact for all $p$ in $U$. Then the
result follows from part (3) of the Lemma, and the five lemma.
\endproof

\subsection{Algebraization}

To show that $\TSM$ is an algebraic stack, we need to
verify  conditions \ref{Sec:stack-criteria}-(\ref{It:formal-compatible}) and
(\ref{It:pro-rep}). 

Condition (\ref{It:pro-rep}) calls for Schlessinger's conditions
as in
\cite{Artin} section 2, (S 1) and (S 2). 

Condition (S 1), and even the stronger
condition (S 1') (see \cite{Artin} (2.3)) follows from standard
principles. Indeed, given an infinitesimal extension $A'\to A \to A_0$, and a
ring homomorphism  $B \to A$ such that $B \to A_0$ is surjective, and given
stable twisted objects $(\xi_{A'}, C_{A'} \to \Spec A',\cA_{A'}) $ and
$(\xi_{B}, 
C_{B} \to \Spec B,\cA_{B}),$ there is a gluing $C \to \Spec A'\times_A
B$, along with gluings of the schemes $U_\alpha$ underlying the charts. The
objects $\eta_\alpha$ can be glued since the condition applies to $\cM$.

 The finiteness condition (S 2)
 follows from the properness of $\cC$ and the exactness of the
pushforward along $\cC \to C$, because of the cohomological description
of ${\rm D}(I)$.

We now come to condition (\ref{It:formal-compatible}) - that formal
deformations 
are algebraizable.

\begin{proposition} Let $A$ be a complete noetherian local ring
over our base scheme $S$, with maximal ideal $\m$. Set $A_n =
A/\m^{n+1}$. Let $\xi_n\colon  {\cC}_n \to {\cM} $ be a
twisted stable  map over $\Spec A_n$, together with
isomorphisms $\alpha_n\colon  \xi_n\rest{\Spec A_{n-1}}
\simeq \xi_{n-1}$. Then there exists a  twisted stable
map
$\xi\colon  {\cC} \to {\cM}$ over
$\Spec A$ together with isomorphisms $\xi\rest{\Spec A_n}
\simeq \xi_n$, compatible with the $\alpha_n$.
\end{proposition}

\proof

{\sc Step 1: construction of $C \to \bM$.} Let $\pi_n\colon  {\cC}_n\to C_n$ be
the moduli 
space of
${\cC}_n$. The existence of the stack $\SM$ insures the existence of a
projective curve $C$ over $\Spec A$, together with
isomorphisms $C\rest {\Spec A_n} \simeq C_n$
compatible with the isomorphisms $C_n\rest{\Spec
A_{n-1}} \simeq C_{n-1}$, and a map $C \to {\bM}$ extending the maps
$C_n \to {\bM}$. 

{\sc Step 2: construction of $\cC \to \cM$.} The stack ${\cM}\times_{{\bM}}C$
is proper and 
tame over $\Spec A$, and its restriction to $\Spec A_n$ is
${\cM}\times_{{\bf M}}C_n$. Consider the maps
$$
\eta_n = \xi_n\times \pi_n\colon {\cC}_n \to {\cM}\times_{{\bf
M}}C_n.
$$
The morphisms $\eta_n$ are representable and finite, so they define a
sheaf of finite algebras ${\eta_n}_*{\cO}_{{\cC}_n}$
over ${\cM}\times_{{\bM}}C_n$. By the existence
theorem \ref{Th:groth-ex} given in the appendix, there is a sheaf ${\cR}$ of
finite  ${\cO}_{{\cM}\times_{{\bM}}C}$-algebras with compatible isomorphisms
${\cO}_ {\cC} \rest {\Spec A_n} \simeq {\cO}_{{\cC}_n}$; by \cite{L-MB}
Proposition~14.2.4, there is a finite representable
map $\eta\colon  {\cC} \to {\cM}\times_{{\bM}}C$
such that $\eta_*{\cO}_{\cC} = {\cA}$.

 By the local criterion of flatness
(\cite{Matsumura}, Theorem~49, condition~5)
${\cC}$ is flat over $A$. By deformation
theory, ${\cC}$ is a nodal stack over $A$.

Let us show that the coarse moduli scheme of ${\cC}$ is
indeed
$C$. Call $\overline C$ the coarse moduli space
of ${\cC}$, with the induced map $\overline C \to
C$. This is finite. Because of the local criterion of
flatness, it is flat (see \cite{Matsumura}, 20.G). Furthermore
since the map $\overline C\rest{\Spec A_0} \to C_0$ is an
isomorphism autside the special locus of $C_0$, the map
$\overline C \to C$ is finite, flat and of degree 1, so it is
an isomorphism.

{\sc Step 3: construction of $\Sigma^{\cC}$.}  Let  $\Sigma^{\cC_n} \subseteq
{\cC}_n$ be 
the union of the 
markings. Again by the existence theorem, there is a unique
substack $\Sigma^{\cC} \subseteq {\cC}$ whose intersection with 
the ${\cC}_n$ is $\Sigma^{\cC_n}$. The stack $\Sigma^{\cC}$ is flat
over $\Spec A$, and is \'etale over the closed point, so it
is \'etale. By the same argument, the diagonal is also
\'etale. Therefore $\Sigma^{\cC}$ is a gerbe over its moduli space,
which is a union of disjoint sections of $C$.

{\sc Step 4: structure of ${\cC} \to C$. }
We only have left to prove that the morphism ${\cC} \to C$
is an isomorphism outside the union  of $\Sigma^{\cC}$ and the singular locus
${\cC}\sing$ of the map ${\cC} \to \Spec A$. There
exists a closed substack of ${\cC}$ where the inertia
groups are nontrivial; this substack meets ${\cC}_0$ at
most at points of $\Sigma^{\cC}$ and ${\cC}\sing$. We need
to analyze the structure of ${\cC}$ near a point of
$\Sigma_0$ and near a point of ${\cC}_0 \cap {\cC}\sing$.

Take a
geometric point $p = \Spec \Omega \to  {\cC}_0$, and let
$R$ the completion of the strict henselization of ${\cC}$
at $p$. Call $\Gamma$ the automorphism group of $p$; this is
a cyclic group of order prime to the characteristic of $A_0$.
Denote  $R_n
= R\otimes_A A_n$.

{\em Case 1:} $p\in \Sigma$.

Now
$R_0$ is isomorphic to $A_0[[\xi_0]]$, where $\xi_0$
is an indeterminate which is a semi-invariant for $\Gamma$.
The ideal generated by $\xi_0$ defines $\Sigma$ at $p$.
Lift $\xi_0$ to a semi-invariant element $\xi\in R$. Then $R =
A[[ \xi]]$; we only need to check that the ideal generated by
$\xi$ defines $\Sigma$. Denote the reduction of $\xi$ in
$R_n$ by $\xi_n$. The ideal of
$\Sigma_n$ in $R_n$ is generated by a semi-invariant
$\zeta_n$, and it is easy to see that $\zeta_n = u \xi_n$,
where $u\in R_n$ is an invariant unit. So $\xi_n$ also
generates the ideal of $\Sigma_n$ in $R_n$, and this implies
that $\xi$ generates the ideal of $\Sigma$ in $R$. This
proves that ${\cC}$ is isomorphic to $C$ outside
$\Sigma$ in a neighborhood of each point of $\Sigma\cap{\cC}_0$.

{\em Case 2:} $p\in  {\cC}\sing$. 

This
time $R_0 $ is isomorphic to $A_0[[\xi_0,
\eta_0]]/(\xi_0\cdot\eta_0)$. We can choose $\xi_0$ and
$\eta_0$ semi-invariants. By deformation theory, there exists
a lifting $\xi'$ of $\xi_0$ and $\eta'$ of $\eta_0$, such that
$\xi'\cdot\eta'\in A$, in other words, $R =
A[[\xi',\eta']]/(\xi'\cdot
\eta' - a)$. Let $\xi$ be any semi-invariant lifting of
$\xi_0$. Then $\xi = u\cdot \xi'$ for some unit $u\in R$. If we
denote $\eta'' = u^{-1}\eta'$, then we have $R =
A[[\xi,\eta'']]/(\xi\cdot
\eta'' - a)$. Denote by $\chi_\xi$ and $\chi_\eta$ the
characters by which  $\Gamma$ acts on $\xi_0$ and $\eta_0$.
Write
$$
\eta = {\frac{1}{|\Gamma|}}\sum_{\gamma\in \Gamma}
\chi_\eta^{-1}(\gamma) \cdot\gamma\eta''.
$$
Clearly $\eta$ is a semi-invariant lifting of $\eta_0$. If
$a=0$, then clearly $\xi\cdot \eta = 0 $ and we are done. If,
on the other hand, $a\neq 0$, then $\chi_\eta =
\chi_\xi^{-1}$, and we have
$$
a = {\frac{1}{|\Gamma|}}\sum_{\gamma\in
\Gamma} \gamma a = {\frac{1}{|\Gamma|}}\sum_{\gamma\in
\Gamma} \gamma(\xi\eta'') = {\frac{1}{|\Gamma|}}\sum_{\gamma\in
\Gamma} \chi_\xi(\gamma) \xi
\cdot\gamma\eta'' =\xi {\frac{1}
{|\Gamma|}}\sum_{\gamma\in
\Gamma} \chi_\eta^{-1}(\gamma)
\cdot\gamma\eta'' = \xi\eta.
$$ 
Again, this
proves that ${\cC}$ is isomorphic to $C$ outside
${\cC}\sing$ in a neighborhood of each point of ${\cC}\sing\cap {\cC}_0$.

 So the pair
${\cC}$, $\Sigma$ with the map ${\cC} \to {\cM} $
gives a twisted stable map. This map is stable, since $C \to
{\bM}$ is stable.\endproof

Thus $\TSM$ is an algebraic stack. Since an automorphism of a twisted object
$(\xi, C, \cA)$ fixing $C\to \bM$ is determined by its action on the generic
object $\xi$, and since $\cM$ is a tame Deligne--Mumford stack, we obtain that
$\TSM 
\to \SM$ is of Deligne--Mumford type and tame.

\begin{remark}
The proof of case (2) above implies the following lemma, which will be used
later (\ref{Prp:balanced-is-clopen}).
\end{remark}
\begin{lemma}\label{Lem:balanced-is-clopen} Let $(\xi, C \to S,\cA)$ be a
twisted object. Then there is an open and closed subscheme $T
\subseteq S$ such that if $s_0$ is a geometric point of $S$,
then the pullback of $(\xi, C \to S,\cA)$ to $s_0$ is balanced
if and only if $s_0$ is in $T$.
\end{lemma}

\section{The weak valuative criterion} \label{Sec:weak-valuative-criterion}
We wish to show that
 $\TSM$ is proper and $\tsm$ is
 projective. 
We start by verifying the weak valuative criterion for
$\TSM$. 

Let $R$ be a discrete valuation ring, $S=\Spec(R)$. Let $\eta\in S$ be the
generic point, $s\in S$ the special point. For a finite extension of
discrete valuation rings
$R\subset R_1$, we denote by $S_1, \eta_1, s_1$ the corresponding schemes.

\begin{proposition} Let $(\xi,C_\eta\to \{\eta\},\cA_\eta)$ be a  stable
twisted 
object. Then there is a finite extension of discrete valuation rings 
$R\subset R_1$ and an extension
$$
\begin{array}{ccc}
(\xi\times_SS_1, C_\eta\times_SS_1,\cA_\eta\times_SS_1)&\subset
&(\xi_1,C_1,\cA_1) 
\\ 
\down		 &	 & \down \\
\{\eta_1\} & \subset & S_1, \end{array}
$$
that is, $(\xi_1,C_1\to S_1,\cA_1)$ is a stable twisted object and its pullback
to $\{\eta_1\}$ is isomorphic to the pullback of $(\xi,C_\eta\to
\{\eta\},\cA_\eta)$ to $\{\eta_1\}$. The
extension is unique up to a unique isomorphism, and its 
formation commutes with further finite extensions of discrete valuation
rings. If  $(\xi,C_\eta\to \{\eta\},\cA_\eta)$ is balanced then $(\xi_1,C_1\to
S_1,\cA_1)$ is 
balanced as well.
\end{proposition}

\proof We will proceed in steps.

{\sc Step 1: extension of $C$.} Let  $f_\eta\colon C_\eta\to \bM$ be the
coarse moduli morphism. By \cite{Fulton-Pandharipande} and
\cite{Behrend-Manin}, we know that $\SM$ is a proper  
algebraic stack. By the weak  
valuative criterion for $\SM$, it 
follows that there 
is a finite extension of discrete valuation rings $R\subset R_1$ and an
 extension 
\begin{equation}\label{Eq:extension-coarse-stable-map}
\begin{array}{ccccc}
C_\eta\times_SS_1 & \subset& C_1 & \stackrel{f_1}{\to} & \bM \\
\down		 &	 & \down & &\\
\{\eta_1\} & \subset & S_1 & &\end{array}
\end{equation}
such that $f_1\colon C_1\to \bM$ is a family of Kontsevich stable maps. The
extension (\ref{Eq:extension-coarse-stable-map}) is unique up to a unique
isomorphism and commutes with further base changes.

We now replace $R$ by $R_1$. Let $\pi\in R$ be a uniformizer.

{\sc Step 2:  extension of $C_\eta\to \cM $ over $C\gen$.} We may assume that
over
$(C_\eta)\gen$ we have a chart with $U_\eta = (C_\eta)\gen$ and
trivial $\Gamma$: such a chart is compatible with any other chart, so we
may add it to the atlas. We may extend this via $U=C\gen$. However the map
$C_\eta\to \cM $ does not necessarily extend over $C\gen$, so we may need a
further base change.

Let $C_i$ be the components of the special fiber $C_s$, and let $\zeta_i\in
C_i$ be
the generic points. Consider the localization $C_{\zeta_i}$. This is the
spectrum of a discrete valuation ring. By the weak valuative criterion for
$\cM$, there is a finite cover $\tilde{C}_i\to C_{\zeta_i}$ and a map
$\tilde{C}_i\to \cM$ lifting the map on $(C_\eta)\gen$.

We proceed to simplify these schemes $\tilde{C}_i$. Denote $S_n = \Spec
R[\pi^{1/n}]$. By Abhyankar's lemma
(\cite{G-SGA1}, exp. XIII section 5) we may assume that $\tilde{C}_i$ is
\'etale 
over $C_{\zeta_i}\times_SS_{n_i}$. By descent we already have a lifting
$C_{\zeta_i}\times_SS_{n_i} \to \cM$. Taking $n$ divisible by all the $n_i$,
we have an extension $C_{\zeta_i}\times_SS_{n} \to \cM$. 

We now replace $R$ by $R[\pi^{1/n}]$. Thus there is an extension $C_{\zeta_i}
\to \cM$. Note that this extension is unique, since $\cM$ is separated. For
the same reason it commutes with further base changes.

By \cite{dJ-O}, there is a maximal open set $U\subset C\gen$ with an extension
of the morphism,
$U \to \cM$. Since $U$ contains both $(C_\eta)\gen$ and $\zeta_i$, we have
$U=C\gen\setminus P$ for a finite set of closed points $P$. By the purity
lemma (Lemma \ref{Lem:purity-lemma})
we have that $U = C\gen$, and there is an extension $C\gen \to \cM$.

The uniqueness of the lifting in the purity lemma guarantees that this
extension is unique up to a unique isomorphism, and commutes with further base
changes. 

{\sc Step 3: extension of $C_\eta\to \cM$ over non-generic nodes.} Let
$p\in \Sing(C)$ be a node, and assume $p$ is not in the closure of
$\Sing(C_\eta)$.

To build up a chart near $p$, we first choose an \'etale neighborhood
$ W\to C$ of $p$, as follows. Let $\bar W$ be a Zariski neighborhood such
that  
$\bar W\cap \Sing(C) = \{p\}$. We already have a morphism $\bar W\to
\cM$. We take 
an \'etale neighborhood
$W$ of $p$ over $\bar W$ so that $p$ is a split node on $W$. Thus
we can find 
elements $s_1,s_2$ in the maximal ideal of $p$, such that $\m_p=(s_1,
s_2, \pi)$, satisfying the equation 
$s_1s_2=\pi^r$ for some $r>1$. In other words, $W$ is also an \'etale
neighborhood of the closed point $\{s_1=s_2=\pi=0\}$ in $W' = \Spec
R[s_1,s_2]/(s_1s_2-\pi^r)$.

Write $r = r_1 r_2$, where $r_2$ is the maximal power of the residue
characteristic in $r$, and thus $r_1$ is prime to the characteristic.

We now find a nodal curve $U_0\to S$, which is Galois over $W$, with an
equivariant extension $U_0 \to \cM$. Let $U_0'\to W'$ be defined as follows:
$U_0'= \Spec R[t_1,t_2]/(t_1t_2-\pi^{r_2})$, and $s_1 = t_1^{r_1}, s_2 =
t_2^{r_2}$. 
There is an obvious {\em balanced} action of the group of $r_1$-th roots of
unity 
$\mmu_{r_1}$ on $U_0'$ via $ (t_1,t_2) \mapsto (\zeta t_1, \zeta^{-1} t_2)$,
and 
$U_0'\to W'$ is the 
associated quotient morphism. Let $U_0 = U_0'\times_{W'} W$. The action of
$\mmu_{r_1}$ clearly lifts to $U_0$. We write $q_0\colon U_0 \to W$ for the
quotient 
map.

Notice that $U_0'$, and therefore $U_0$, satisfies Serre's condition
$S_2$. Also, 
$U_0'$ (and therefore also $U_0$) has a homeomorphic nonsingular cover $\tilde
U'$ given by taking  roots of order $r_2$ of $t_1$ and $t_2$. Thus the local
fundamental group is trivial, and the purity
lemma applies. By the purity
lemma we have a lifting $U_0\to \cM$ of the morphism $f\circ q_0 \colon  U_0
\to 
\cM$. Notice that for $g\in \mmu_{r_1}$ we have $f\circ q_0 = f\circ q_0 \circ
g$. Thus by the uniqueness of the lifting, the map $U_0 \to \cM$ commutes
with the action of $\mmu_{r_1}$.

Let $\eta_0\in \cM(U_0)$ be the associated object and
let 
$\eta_p$ be the fiber over $p$. The group $\mmu_{r_1}$ acts on $\eta_0$,
stabilizing $\eta_p$. Let $\mmu_{r_0} \subset \mmu_{r_1}$ be the subgroup
acting 
trivially on $\eta_p$. Denote by $U = U_0/\mmu_{r_0}$ and  $\Gamma=
\mmu_{r_1}/\mmu_{r_0}$ the quotients. By 
lemma \ref{Lem:descent}  there is a $\Gamma$-equivariant object $\eta\colon U
\to \cM$, and the action is essential. 

Thus we have a chart $(U, \eta, \Gamma)$ near $p$. This is easily seen to
be unique and to commute with further base changes. Note that the chart we have
constructed is automatically balanced.

{\sc Step 4: extension of charts over generic nodes.} Let
$p_\eta\in \Sing (C_\eta)$, and let $p\in C_s$ be in the closure of $p_\eta$.
The construction of a chart here is similar to step 2. We choose a Zariski
neighborhood $\bar W$ of $p$ such that $\bar W\cap \Sing(C) = \{p_\eta,p\}$. As
before, we may choose an \'etale neighborhood $W$ over $\bar W$ such that
$p_\eta$ is a split node, thus $W$ is \'etale over $\{s_1s_2=0\}$. We may
assume that $U_\eta$ is given by $t_1t_2=0,$ where $t_i^r = s_i$ for an
appropriate integer $r$ prime to the residue characteristic, with the group
$\Gamma = \mmu_r$. Otherwise we can 
add a compatible chart with such a $U_\eta$ to our atlas. There is an obvious
extension $U_\eta \subset U$ via the same equation $\{t_1t_2=0\}$, and the
action of $\Gamma$ extends automatically. We
already have a
lifting $U\setminus \{p\} \to \cM$. Since $U$ is normal crossings, the purity
lemma applies, and guarantees that there is a unique equivariant extension
to $U$. It is easy to check that this gives an extension of the chart, which is
essential (as the $\Isom$ scheme of $\cM$ is finite unramified). Also, if the
chart on $(\xi, C_\eta \to S,\cA_\eta)$ is balanced at $p_\eta$ then the chart
we have 
constructed is evidently balanced as well.

{\sc Step 5: extension of $C_\eta\to \cM$ over $\Sigma_i$.} This step
is identical to the previous one.

The uniqueness and base change properties are straightforward. \endproof

\section{Boundedness}

\subsection{The statement} The bulk of  this section is devoted to proving the
following result: 

\begin{theorem}\label{Th:boundedness}  Given a morphism $T \to \SM$, where $T$
is a 
scheme of finite type over $\kappa$, there exists a morphism
$T' \to T$ of finite type, and a lifting $T' \to
\TSM$, such that a geometric point of $\TSM$
is in the image of  
$T'$ if its image in $\SM$ is in the image of $T$.
\end{theorem}

Let $C \to T$ be the corresponding underlying family of curves, with stable map
$f \colon  C \to {\bM}$. 

\subsection{The smooth case}
\begin{proposition} The Theorem holds when $C \to T$ is
smooth.
\end{proposition}

By noetherian induction, it suffices to prove the Proposition
after a dominant base change of finite type on $T$. In
particular we may assume $T$ irreducible. Consider $\overline
{f(C)}$, the closure of the image of $C$ in ${\bM}$; there
exists a finite map $W \to (\overline {f(C)} \times_ {\bM}{\cM})_\red$
which is generically \'etale over $(\overline 
{f(C)} \times_ {\bf M}{\cM})_\red$ (\cite{L-MB}, 16.6.). Then by
Lemma \ref{Lem:tame-cms-pullback}, the scheme $\overline {f(C)}$ is the moduli
space of 
$(\overline {f(C)} 
\times_ {\bf
M}{\cM})_\red$, so, by \cite{L-MB} 11.5, $(\overline {f(C)} \times_
{\bM}{\cM})_\red$ is generically \'etale over $\overline
{f(C)}$. So it follows that $W$ is also generically \'etale
over $\overline {f(C)}$.

Let $F \to W \times_{{\bM}} C$ be the normalization of an
irreducible component of $W \times_{{\bM}} C$ mapping
surjective onto $C$. By refining $T$, we may assume that
\begin{enumerate}
\item  $F \to T$ is smooth;

\item  there exists an open subset $C_0 \subseteq C$ such
that
\begin{enumerate}
\item  $C \setminus C_0$ is a union of disjoint sections
$\Xi$ containing $\Sigma$,

\item  if we denote by $h \colon  F \to C$, then $F_0 \eqdef
h ^{-1}(C_0) \to C_0$ is \'etale, and

\item  the scheme of automorphisms of $\beta\colon  F \to
{\cM}$ obtained by composing $F \to W \to {\cM}$  is
\'etale over $F_0$ of constant degree $\alpha$.

\end{enumerate}
\end{enumerate}

\begin{lemma} \label{Lem:exhausting-covers} 
There exists a morphism $T_1 \to T$ of
finite type, and a morphism $E \to F_{T_1}$ such that $E \to
T_1$ is smooth of relative dimension 1, $E \to F_{T_1}$ is
\'etale over $(F_0)_{T_1}$ of degree $\alpha$, and tame over
$(F \setmin F_0)_{T_1}$, with the following property.

For every geometric point $t \colon   \Spec \Omega \to T$, and
every tame cover $E' \to F_t$ of order $\alpha$ which is
\'etale over $(F_0)_t$, there exists a lifting $t_1 \colon 
\Spec\Omega \to T_1$, such that $E' \to F_t$ is isomorphic to
$E_{t_1}$ as a covering of $F_t$.
\end{lemma}

\proof Since $E' \to F_t$ is tame of degree $\alpha$, the
genus of every connected component
of $E'$ is bounded by some integer $g_e$. It is sufficient
to consider one connected component of $E'$ at a time; call
$\alpha'$ its degree over $F_t$. Consider the stack of maps
${\cK}_{g_e}(F/T, [\alpha'])$; there is a locally closed
substack ${\cK}_0$ consisting of maps which are \'etale
over $F_0$ and with smooth source curve. 
There exists a scheme $T' \to {\cK}_0$
surjective and of finite type; let $\overline E \to F_{T'}$ be
the corresponding universal family. We use the following well-known lemma:

\begin{lemma}\label{Lem:tame-constructible} There exists a constructible subset
$T^{\rm 
tame} \subseteq T'$ such that for a geometric point
$t'\colon   \Spec\Omega \to T'$ the fiber $\overline E_{t'} \to
F_{t'}$ is a tame covering  if and only if $t'$ is in $T^{\rm
tame}$.
\end{lemma}

\proof 
By Noetherian induction it suffices to show that there is an open set where the
lemma holds, and thus we may replace $T'$ by a dominant generically finite
scheme over it. Let $\bar\eta$ be the generic point of $T'$, and let $\tilde
E_{\bar\eta} 
\to F_{\bar\eta}$ be the normalized ${\mathfrak S}_n$ cover associated with
$E_{\bar\eta} \to F_{\bar\eta}$, where ${\mathfrak S}_n$ is the symmetric
group. There is a finite extension $\eta'' \to \eta$ such that the ${\mathfrak
S}_n$ cover $\tilde E_{\bar\eta} \to F_{\bar\eta}$ comes from   $\tilde
E_{\eta''} \to 
F_{\eta''}$. There is a dominant generically finite $T'' \to T'$ such that
$\eta_1$ is the generic point of $T''$ and $\tilde E_{\eta''} \to
F_{\eta''}$ comes from   $\tilde E_{T''} \to
F_{T''}$ and the latter is the normalized  ${\mathfrak S}_n$
cover associated with $E_{T''} \to
F_{T''}$. The set of fixed points of elements of order $p$ in ${\mathfrak S}_n$
is closed in $\tilde E_{T''}$ , and the image in $T'$ is
constructible. \endproof
 
To conclude the proof of Lemma \ref{Lem:exhausting-covers}, we take $S_1$ to be
a disjoint union of locally closed subschemes of $T'$ covering
$T^{\rm tame}$. \endproof

Consider $\widetilde E_2 \eqdef (E \times_{C_{S_1}} E)\norm$,
$\widetilde E_3 \eqdef (E \times_{C_{S_1}} E \times_{C_{S_1}}
E)\norm$. By Abhyankar's Lemma, $\widetilde E_2,
\widetilde E_3 \to T_1$ are smooth.

We have two morphisms $\alpha_1, \alpha_2\colon  \widetilde
E_2 \to {\cM}$ obtained by composing the map $E \to F \to W
\to {\cM}$ with the two projections $\widetilde E_2 \to
E$. Consider the scheme $\widetilde I_2 = \Isom_{\widetilde
E_2}(\alpha_1, \alpha_2)$. The morphisms $\widetilde I_2 \to
\widetilde E_2$ is finite, and therefore $\widetilde I_2$ is
projective over $T_1$. So there exists a scheme $T_2 =
\Sect_{T_1}( \widetilde I_2 \to \widetilde E_2)$, quasi-projective over
$T_1$, parametrizing sections of $\widetilde I_2 \to
\widetilde E_2$. In other words, there exists a diagram
$$
\begin{array}{c}\widetilde I_2\rest {T_2}\\
 \downarrow \curveuparrow \\
\widetilde E_2\rest {T_2}\\ \downarrow\\ T_2;
\end{array}
$$
denote the universal section $\sigma\colon  \widetilde E_2\rest
{T_2} \to \widetilde I_2\rest {T_2}$.

Considering the three natural projections $p_{12}, p_{13},
p_{23}\colon  \widetilde E_3 \to \widetilde E_2$, we have a
composition map $\mu\colon  p_{12}^*\widetilde I_2 \times
p_{23} ^*\widetilde I_2 \to p_{13} ^*\widetilde I_2$. There
exists a closed subscheme $T_3 \subseteq T_2$ where the
composition
$$
\mu\circ(p_{12}^*\sigma \times p_{23}^*\sigma)
\colon  \widetilde E_3 \to p_{12}^*\widetilde I_2 \times
p_{23} ^*\widetilde I_2 \to p_{13} ^*\widetilde I_2\quad
\hbox{equals}
\quad p_{13} ^* \sigma \colon   \widetilde E_3 \to p_{13}
^*\widetilde I_2.\leqno (*)
$$

\begin{lemma} There exists  a lifting $\lambda\colon 
(C_0)_{T_3} \to {\cM}$ of $C \to {\bM}$, such that the following holds: 
whenever $t \colon   \Spec\Omega \to T$ is a geometric point,
and $\tau\colon  (C_0)_t \to{\cM}$ is a lifting of $C_t\to
{\bM}$, there exists a lifting $t_3 \colon  \Spec\Omega \to
T_3$ and an isomorphism $t_3 ^*\lambda \simeq \tau$.
\end{lemma}

\proof The section $\sigma_3 \eqdef \sigma\rest{T_3}\colon 
(\widetilde E_2)_{T_3} \to (\widetilde I_2)_{T_3}$, when
restricted to the open set $E_0:=$ the inverse image of $C_0$,
gives  descent data for constructing $\lambda\colon 
(C_0)_{T_3} \to {\cM}$, which is effective since $\cM$ is a stack. Given a
geometric point $t \colon  
\Spec\Omega \to T$ and a lifting $\tau\colon  (C_0)_t \to{\cM}$ of
$C_t\to {\bM}$, we have a diagram 
$$
\begin{array}{ccccc}(F_0)_t & \subseteq & F & \to & W \\
&&&&\down\\
\down&&\down&& {\cM}\\
&&&&\down\\
(C_0)_t & \subseteq & C & \to & {\bM} \end{array}
$$
giving two morphisms $\beta, \beta_2\colon  (F_0)_t \to {\cM}$ which
induce the same morphism $(F_0)_t \to {\bM}$. 
Therefore the scheme $\Isom_{(F_0)_t}(\beta, \beta_2)$ is
\'etale over $F_0$, as it is a torsor under $\Aut_{(F_0)_t}(
\beta)$. For the same reason, it is tame along $(F \setmin
F_0)_t$. Thus there exists a lifting $t_1 \colon  \Spec\Omega
\to T_1$ such that $E_{t_1} \to F_{t_1}$ is isomorphic to the
normalization of $F_{t_1}$ in $\Isom_{(F_0)_t}(\beta,
\beta_2)$. In particular, we have an isomorphism of
$\beta\rest{E_{T_1}}$ with $\beta_2$; by conjugation this
gives a section $\left((\widetilde E_2)_{t_1}\right)\rest{E_0}
\to \widetilde I_2$, which extends to $(\widetilde
E_2)_{t_1} \to \widetilde I_2$, since $(\widetilde
E_2)_{t_1}$ is a smooth curve. This exended morphism is
easily seen to satisfy the cocycle condition $(*)$, so it
gives a point $t_3 \colon   \Spec\Omega \to T_3$, and the
corresponding descent data  give the lifting $\tau\colon 
(C_0)_t \to {\cM}$.

Replacing $T_3$ by a regular stratification, this proves the following:

\begin{lemma} There exists a  morphism $T_0 \to T$
of finite type, with $T_0$ regular, and a lifting $(C_0)_{T_0}
\to {\cM}$, such that whenever $t \colon   \Spec\Omega \to
T$ is a geometric point, and $\tau\colon  (C_0)_t \to{\cM}$
is a lifting of $C_t\to {\bM}$, there exists a lifting $t_0
\colon  \Spec\Omega \to T_0$ and an isomorphism $t_0 ^*\lambda
\simeq \tau$.
\end{lemma}

\begin{lemma} Let $C \to \Spec \kappa$ be a smooth
curve over an algebraically closed field, $C_0 = C \setmin
\{ p_1 ,\ldots ,p_n\}$ an open subset, $C_0 \to {\cM}$ a
morphism. Then there exists a twisted pointed curve $( {\cC},
\overline \Xi_1,\ldots , \overline \Xi_n)$ with 
moduli space $(C,p_1 ,\ldots ,p_n)$, and a representable
morphism $ {\cC} \to  {\cM}$.
\end{lemma}

\proof Choose an \'etale neighborhood $U_i \to {\cM}$ of
the inverse image in $ {\cM}$ of $f(p_i)$, such that $U_i
\times_{\bM} {\cM}= [V_i / \Gamma_i]$, where $\Gamma_i$
is a finite group. Consider
$$
V_i \times_{\bM} C \to U_i \times_{\bM} C \to C.$$

Passing to an \'etale neighborhood $U_i' \to U_i \times_{\bf
M} C \to C$ of $p_i$, we may assume that the normalization
$V_i'$ of some irreducible component of $U_i' \times_{U_i}
V_i$ has exactly  one point $q_i$ over $p_i$. Let $\Gamma_i' $
be the stabilizer of $V_i'$. Note that $\Gamma_i' $ fixes
$q_i$.  Consider the subgroup
$\Gamma_i'' \subset \Gamma_i'$ stabilizing the object $q_i
\to {\cM}$. Denote $ \overline V_i' = V_i' / \Gamma_i''$, $
\overline \Gamma_i' = \Gamma_i'  / \Gamma_i''$, and $
\overline q_i $ the image of
$q_i$. Refining
$U_i'$ if necessary, we may assume that $\overline q_i $ is
the unique fixed point of $\overline \Gamma_i' $ in
$\overline V_i'$. By Corollary \ref{Cor:descent-to-quotient} we have that
$$
(\overline V_i'  \to C, \overline V_i'  \to \cM, \overline \Gamma_i' ) 
$$
forms a chart at $p_i$. \endproof

Back to the proposition.
Let $ \overline \eta \to T_0$ be the geometric  generic
point. By the lemma, there exists a representable lifting
$$
{\cC}_{ \overline \eta} \to  {\cM}
$$
over a twisted curve $({\cC}_{ \overline \eta},
({\Xi_1})^{\cC_{ \overline \eta}},\ldots,({\Xi_n})^{\cC_{ \overline
\eta}})$.

There exists a dominant quasifinite morphism $T_1 \to T_0$
such that the lifting above descends to $T_1$. We claim that
there exists an open dense subscheme $T_2 \subset T_1$ such
that
$ {\cC}_{T_2}$ is representable. This follows from the
following well known fact:

\begin{lemma} Let
$f\colon  {\cX}_1  \to  {\cX}_2$
be a morphism of Deligne--Mumford stacks, and let $ {\cX}\to
S$ be proper. Then there is an open subset $S_0 \subset S$
such that if $s\colon \Spec \Omega \to S$ is a geometric point, the
restricted morphism $f_s\colon  ({\cX}_1)_s \to {\cX}_1$
is representable if and only if $s$ factors through $S_0$.
\end{lemma}

\proof Consider the morphism of inertia stacks $ {\cI}_{
{\cX}_1} \to f^* {\cI}_{
{\cX}_2}$. The kernel is a finite unramified representable
group stack $ {\cK} \to {\cX}_1$, therefore the
complement $ {\cK}^0 \subset {\cK}$ of the identity is
also finite over $ {\cX}_1$. Then $S_0$ is the complement
of the image of $ {\cK}_0$ in $S$, which is closed.
\endproof

We may now assume that we have a lifting $ {\cC} \to
{\cM}$ over $T_1$. To conclude the proof of the
proposition, there exists an open and closed subscheme $T'
\subset T_1$, such that for any section $\Xi_i$ in $\Xi
\setmin \Sigma$ the inertia group of a geometric point in
$\Xi_i^\cC$ is trivial if and only if the point lies over $T'$.
\endproof

\subsection{Proof of the Theorem}
 By taking a surjective map of
finite type $T_1 \to T$ we may assume that
\begin{enumerate}
\item  $C \to T$ is of locally constant topological type,

\item  All irreducible components of $C$ have
geometrically irreducible fibers,

\item  Any irreducible component $N_i$ of the singular
locus of
$C \to T$ maps isomorphically to its image.

\end{enumerate}

This implies that the normalization $C\norm \to C$ can be
viewed as a union of families of smooth pointed curves with a
stable map to $ {\bM}$. By the proposition we may assume
that there is a dominant morphism $T_0 \to T$
of finite type, and a lifting
$$
\begin{array}{ccc}
  {\cD}           & \stackrel{F}{\lrar} & {\cM}  \\
   \down          &     &  \down \\
   (C_{T_0})\norm & \to &  {\bM} \\
   \down          &     &        \\
   T_0.           &     &
\end{array}
$$

For each $N_i$, denote $ \overline \Sigma_{i,1}, \overline
\Sigma_{i,2} \subset C\norm$ be the two sections over $N_i$,
and $  \Sigma_{i,1},
\Sigma_{i,2} \subset {\cD}$ the gerbes in $ {\cD}$
over them. replacing $T_0$ by an \'etale cover we may assume
that $  \Sigma_{i,1} \to T_0$ and $  \Sigma_{i,2} \to T_0$
have sections $\sigma_{i,1}\colon  T_0 \to  \Sigma_{i,1}$ and
$\sigma_{i,2}\colon  T_0 \to  \Sigma_{i,2}$. Refining $T_0$ we may
assume that there are charts
$$(U_{i,j},\eta_{i,j}, \Gamma_{i,j}) $$
 for $
{\cD}$ along $  \Sigma_{i,1} $ and $  \Sigma_{i,2} $ and
that the sections $\sigma_{i,j} $ lift to
$\widetilde \sigma_{i,j}\colon T_0 \to U_{i,j}$. Also we may assume
that $\widetilde \sigma_{i,j}$ surjects onto the inverse
image of $\sigma_{i,j}$ in $U_{i,j}$. This implies that the
groups $\Gamma_{i,j}$ are cyclic, and isomorphic to $\Aut
(\sigma_{i,j})$.

Denote $\tau_{i,j} = F\circ \sigma_{i,j}\colon  T_0 \to
\Sigma_{i,j} \to  {\cM}$. we have a canonical isomorphism
$\tau_{i,j} = \widetilde \sigma_{i,j}^*\eta_{i,j}$, thus
there is an embedding of $\Gamma_{i,j}\times T_0\subset
\Aut(\tau_{i,j})$.   Consider the subscheme
$T_1$ of the scheme of isomorphisms
$\Isom_{T_0}(\tau_{i,1},\tau_{i,2})$ corresponding to
isomorphisms which identify
$\Gamma_{i,1}$ with $\Gamma_{i,2}$. We may assume that
$T_1$ is a union of disjoint copies of $T_0$. Consider a
section
$\psi\colon  T_0\to T_1$

We choose:
\begin{enumerate}
\item  a geometric point $t_0 \colon  \Spec \Omega \to T_0$,

\item a universal deformation space for $ {\cM}$ at
$t_0^*\tau_{i,1}$, equivariant under the action of
$\Aut(t_0^*\tau_{i,1})$.
\end{enumerate}

By refining $T_0$ and $U_{i,j}$, we may assume that there is
a lifting $\nu_{i,1}\colon  U_{i,1} \to W$ of $\eta_{i,1}$.

Now the chosen isomorphism $\psi$ of $\eta_{i,1}$ with
$\eta_{i,2}$ and the universal property of $W$ imply
that, after refining $U_{i,2}$, there is a lifting
$\nu_{i,2}\colon  U_{i,2} \to W$ of      $\eta_{i,2}$ such that
the two morphisms $\nu_{i,j} \circ \eta_{i,j}\colon  T_0 \to W$
coincide.

We thus obtain a morphism
$$
\eta_i\colon  U_i \eqdef U_{i,1}
\mathbin{\mathop{\cup}\limits_{T_0}} U_{i,2} \to W.
$$

The identification of $\Gamma_{i,1}$ with $\Gamma_{i,2}$
given by $\psi$ defines an action of $\Gamma_{i,1}$  on both
$U_i$. It acts on $W$ via the embedding $\Gamma_{i,1} \subset
\Aut(t_0^*\tau_{i,1})$. The morphism
$\eta$  is
$\Gamma_{i,1}$-equivariant, therefore
$$
(U_i, \eta_i, \Gamma_{i,1})
$$
is a chart for a twisted map $ {\cC}_\psi\to  {\cM}$.

Thaking the union over all $\psi\colon  T_0 \to T_1$,
this defines
$$
\begin{array}{ccc}
 {\cC}          & \to & {\cM} \\
 \down          &     &  \down          \\
 (C_{T_1})\norm & \to &  {\bM} \\
 \down          &     &            \\
 T_1.           &     &        
\end{array}
$$
This  clearly exhausts all liftings over all geometric points. Setting $T' =
T_1$ completes the 
proof of  Theorem \ref{Th:boundedness}.
\endproof
\subsection{Finiteness of fibers of $\tsm\to \sm$}

\begin{lemma}\label{Lem:finite-fibers} Let $\cC$ be a twisted curve over  an
algebraically closed 
field with coarse moduli space $C$. A morphism
$C\to \bM$ can come from at most finitely many representable maps $f:\cC \to
\cM$,  up to isomorphism.
\end{lemma}

\proof First of all observe that given an object $\xi\in \cM(C_\smooth)$  there
are at most finitely many twisted objects on $\cC$ 
whose restriction to $C_\smooth$ is isomorphic to $\xi$; so it is enough to
prove 
that any given map $f\colon T\to \bM$, where $T$ is a smooth curve, comes
from at most finitely many objects on $T$.

Over $\CC$, it is easy to deduce this finiteness result analytically: given a
base point  $t\in T$ and a fixed lifting $\tilde{f}\colon T\to \cM$, it is easy
to associate, in a one to one manner, to  any lifting  $\tilde{f'}\colon
T\to \cM$ a monodromy representation $\phi_{f'}: \pi_1(T) \to \Aut
\xi_t$. Since 
$\pi_1(T)$ is finitely generated there are only finitely many such
representations. 

We now give a short algebraic argument using the language of stacks.

 Call $\cT$ the normalization
of the pullback of $\cM$ to $T$; if there is at least one map $T\to \cM$
giving rise to $f$ then this induces a section $s\colon T\to {\cT}$.
It is proved in \cite{Vistoli:chow-stack} that this section is \'etale; the
argument is similar to that of the purity lemma: by passing to the strict
henselization one can lift $s$ to the universal deformation space. The image of
this lifting is then evidently \'etale, and therefore $s$ is \'etale.
This strongly restricts the structure of $T$. 

 Call $R$ the fiber
product $T\times_{\cT} T$; then the \'etale finite groupoid $R\double T$
is a presentation of $\cT$. Each lifting $T\to \cM$ of $f$ comes
from a section $T\to {\cT}$, so it is enough to show that up to
isomorphisms there are only finitely many such sections. But such a section
$\sigma\colon T\to {\cT}$ is determined by the pullback $U$ of
$T\stackrel{s}{\to} {\cT}$ to $T$ via $\sigma$, and a morphism of groupoids
from $U\times_T U\double U$ to $R\double T$ over $T$. But up to isomorphism
there are only finitely many \'etale coverings $U\to T$ of bounded degree,
because the 
fundamental group of $T$ is topologically finitely generated. For each
$U\to T$ there are only finitely many liftings $U\times_T U\to R$.\endproof

\subsection{Conclusion of proof of the main theorem}
 Theorem
\ref{Th:boundedness} and the fact that the stack $\SM$ is of finite type (see,
e.g., \cite{Behrend-Manin}) imply that the stack $\TSM$ is of finite
type. Since the stack has finite diagonal and satisfies the weak valuative
criterion for properness, it is proper. The induced morphism $\tsm\to \sm$ of
coarse moduli spaces is proper with finite fibers. Since $\sm$ is projective,
we have that $\tsm$ is projective as well. We have already seen that $\TSM \to
\SM$ is of Deligne--Mumford type and tame.
\endproof

\section{Some properties and generalizations}
\subsection{Balanced maps} 
\begin{proposition}\label{Prp:balanced-is-clopen}
The stack $\TSMB$ of balanced twisted stable maps is an open and closed
substack in $\TSM$. 
\end{proposition}

This is immediate form the Lemma \ref{Lem:balanced-is-clopen}.

\subsection{Fixing classes in Chow groups modulo algebraic equivalence} Assume
$\bbS$ is the spectrum of a field. Given a
class $\beta$ in the Chow group of cycles of dimension 1 on $\bM$ modulo
rational equivalence, 
such that its intersection number with the Chern class of the chosen ample
sheaf is $d$, 
we have an open and closed substack $\KO{g}{n}{\bM}{\beta}\subset
\KO{g}{n}{\bM}{d}$. Pulling back, we have an open and closed substack
$\KO{g}{n}{\cM}{\beta}\subset \KO{g}{n}{\cM}{d}$.  A similar construction can
be obtained using classes on $\cM$, whenever one can construct them with
integral coefficients, see, e.g.~\cite{Edidin-Graham}.

\subsection{Stable maps over a base stack}
Suppose we are given a tame morphism $\cM \to \cS$ of algebraic stacks, where
$\cS$ is Noetherian. Suppose further 
that it factors as $\cM \to \bM \to \cS$, where $\cM \to  \bM$ is proper and
quasifinite, and $\bM \to \cS$ is projective (but we do not assume that $\bM$
is representable). We can define a category $\KO{g}{n}{\cM/\cS}{d}$ of
twisted stable maps exactly as before. We claim that the analogue of the main
theorem still holds, which may be of use in some applications. 

The point is the following: most arguments in this paper go through word for
word. To produce a smooth parametrization take a smooth surjective morphism
$\bbS \to 
\cS$, and denote by $\cM_1$ the pullback of $\cM$ to $\bbS$. Then $$
\KO{g}{n}{\cM_1/\bbS}{d} \to  \KO{g}{n}{\cM/\cS}{d}$$ is  representable,
smooth,  surjective  and of finite type. Composing with a smooth
parametrization of the stack $\KO{g}{n}{\cM_1/\bbS}{d}$ we get what we want. To show that
$\KO{g}{n}{\cM/\cS}{d}\to \cS$ is proper, note that its pullback after a
faithfully  flat base change to a scheme is proper.

\subsection{When $\bM$ is only a proper algebraic space}
So far we have assumed that $\bM \to \bbS$ is projective. In case $\bM \to
\bbS$ is only  proper, one can proceed as follows.

First, it was pointed out by Johan de Jong that in this case
$\KO{g}{n}{\bM}{\beta}$ is an algebraic stack, locally of finite type, such
that each irreducible component is proper over $\bbS$. The only part which is
not evident is the properness of the irreducible components. This can be seen
as follows: let $\eta$ be a generic point of such an irreducible component,
with a   stable map $(C_\eta\to \{\eta\}, \Sigma_i, C_\eta \to \bM)$. An easy
gluing arguments allows us to reduce to the case where $C_\eta$ is
irreducible. We then 
replace $\bM$ by the closure of the image of  $C_\eta$. There is a birational
morphism $\bM' \to \bM$ such that $\bM'\to \bbS$ is projective. There is a
lifting $C_\eta \to \bM'$, which has a certain degree $d$ with respect to some
ample sheaf on $\bM'$. The map $C \to \bM'$ induces
$\eta \to\KO{g}{n}{\bM'}{d}$. Take the  closure $\cK$ of the image, which is a 
proper stack. There is a
universal stable map $(C' \to \cK, \Sigma_i', C'\to \bM')$. The composite  map
$(C' \to \cK, \Sigma_i', C'\to \bM' \to \bM)$ may be unstable, but using
Knudsen's  contraction procedure we obtain an associated map $(C \to \cK,
\Sigma_i, C\to\bM)$, giving a morphism $\cK \to  \KO{g}{n}{\bM}{\beta}$ whose
image, which is the component with which we started, is proper. 

Now considering a tame stack $\cM$ proper over $\bbS$, the arguments of the main
theorem along with the one 
above show that $\KO{g}{n}{\cM}{\beta}$  is again  an algebraic stack, locally
of finite type, such 
that each irreducible component is proper over $\bbS$.

\subsection{Other open and closed loci}

Consider a twisted object $(\xi, C \to S,\cA)$, and
take a geometric point $c_0$ of the support $\Sigma^C_i$ of the
$i^{\rm th}$ section $S \to C$. If $(U,\eta, \Gamma)$ is a
chart in the atlas, and $u_0$ is a geometric point of $U$
lying over $c_0$, the stabilizer of $u_0$ in $\Gamma$ has an
order which is independent of the chart. We call this order
the {\em local index} of the twisted object at $c_0$; we can
think of it as a measure of how twisted the object is around
$c_0$. One immediately checks that the local index only
depends on the image of $c_0$ in $S$; this way we get a
function $\epsilon_i\colon  S \to {\bf N}$.

\begin{proposition} The function $\epsilon_i\colon  S \to
{\bf N}$ is locally constant.
\end{proposition}

\proof Let $c_0$ be a geometric point of $\Sigma^C_i$, $(U,\eta,
\Gamma)$ a chart, $u_0$ a geometric
point of $U$ lying over $c_0$. By refining the chart we may
assume that $u_0$ is a fixed point of\/ $\Gamma$. Since the
morphism $\Sigma^C_i \to S$ is \'etale, it follows that $\Gamma$
leaves the whole component of $\Sigma^C_i$ containing $u_0$
fixed, and the thesis follows easily from this.\endproof

\begin{definition} Let ${\mathbf e} = (e_1, \ldots, e_n)$ be a
sequence of positive integers. A twisted object has\/
{\rm global index} equal to ${\mathbf e}$ if the value of the
associated function
$\epsilon_i$ on $S$ is constant equal to $e_i$ for all $i = 1,
\ldots,n$.
\end{definition}

Sometimes the local index can be further refined. For
example, if $G$ is a finite group, and ${\cB}G$ is the
classifying stack of $G$ over  $\bbS$, that is,
the stack whose objects are Galois \'etale covers with
group $G$, then the stabilizer of a geometric point of
$\Sigma^U_i$ is a cyclic subgroup of the automorphism group of an
\'etale Galois cover of the spectrum of an algebraically
closed field. This group is isomorphic to $G$, and the
isomorphism is well defined up to conjugation, so from each
geometric point of $\Sigma^U_i$ one gets a conjugacy class of
cyclic subgroups of $G$, which only depends on the index $i$
and on the image of the geometric point in $S$. It is easy to
check that this conjugacy class is locally constant on $S$.

\section{Functoriality of the stack of twisted stable maps}

\subsection{statements}
Let $\cC \to T$ be a proper $n$-pointed twisted curve of genus $g$ and let $\cC
\to \cM$ be a morphism 
of stacks. Consider the corresponding morphism of coarse moduli spaces $C \to
\bM$. Assume that either $C \to M$ is nonconstant or $2g-2+n>0$. It is well
known (see \cite{Behrend-Manin},  Proposition 3.10) that there exists a 
canonical proper surjective morphism $C \to C'$ of pointed curves which gives a
factorization $C 
\to C' \to \bM$ such that $C' \to M$ is a stable pointed map and the fibers $C
\to C'$ 
are the non-stable trees of rational curves.

\begin{proposition} There exists a  factorization $\cC \to \cC' \to
\cM$ such that 
 	\begin{enumerate}
  		\item $\cC'\to \cM$ is a twisted stable map, and
		\item on the level of coarse moduli spaces this induces the
  		factorization $C \to C' \to \bM$.
	\end{enumerate} The factorization is unique up to a unique isomorphism.
	The formation of $\cC'$ commutes with base change on $T$.
	Moreover, if $\cC \to T$ is balanced, then so is $\cC'$.
\end{proposition}

We denote $\KO{g}{n}{\cM}{*} = \cup_d\KO{g}{n}{\cM}{d}$.
The following corollaries follow immediately from the proposition: 
\begin{corollary}
Let $\phi:\cM \to \cM'$ be a morphism of proper tame stacks over $\bbS$ having
projective coarse moduli spaces. 
Denote by $\cK^{[\phi]}_{g,n}(\cM,d)\subset \KO
{g}{n}{\cM}{d}$ the following  open and
closed locus: if $2g-2+n>0$ then $\cK^{[\phi]}_{g,n}(\cM,d)=
\KO{g}{n}{\cM}{d}$. Otherwise it is the locus where composite map of moduli
spaces $C \to 
\bM \to \bM'$ is non-constant. 
 Then there exists a morphism $\cK^{[\phi]}_{g,n}(\cM,*) \to
\KO{g}{n}{\cM'}{*}$, sending balanced maps to balanced maps, which makes the
following diagram commutative: 
$$ \begin{array}{ccc}
    \cK^{[\phi]}_{g,n}(\cM,*) & \to &\KO{g}{n}{\cM'}{*} \\
      \dar             &     & \dar \\
    \cK^{[\phi]}_{g,n}(\bM,*) & \to &\KO{g}{n}{\bM'}{*}
\end{array} $$
where the vertical arrows are the canonical maps described in Theorem
\ref{Th:stable-maps}
and the morphism   $\cK^{[\phi]}_{g,n}(\bM,*)  \to \KO{g}{n}{\bM'}{*}$ is the
one described by Knudsen and Behrend-Manin, induced by the contraction of
rational trees which become non-stable over $\bM'$.
\end{corollary}

\begin{corollary}
Let $\cK'_{g,n}(\cM,d)\subset \KO{g}{n}{\cM}{d}$ be the open and closed
substack where the last marking $\Sigma_n$ is representable.  Assume either
$d>0$ or $2g-2+n>1$. Then there exists a morphism $\cK'_{g,n}(\cM,d)\to
\KO{g}{n-1}{\cM}{d}$, sending balanced maps to balanced maps, which makes the
following diagram commutative: 
$$ \begin{array}{ccc}
    \cK'_{g,n}(\cM,d) & \to &\KO{g}{n-1}{\cM}{d} \\
      \dar             &     & \dar \\
    \cK'_{g,n}(\bM,d) & \to &\KO{g}{n-1}{\bM}{d}
\end{array} $$
where the vertical arrows are the canonical maps described in Theorem
\ref{Th:stable-maps}
and the morphism   $  \cK'_{g,n}(\bM,d)  \to \KO{g}{n-1}{\bM}{d}$ is the
one described by Knudsen and Behrend-Manin, induced by the contraction of
rational components which become non-stable when forgetting the last marking 
$\Sigma_n^C$. 
\end{corollary}

\begin{remark}
Given a $T$-scheme of finite type $U \to T$, not necessarily proper, one
defines a stable divisorially $n$-marked map $(C \to
U,\Sigma_1,\ldots,\Sigma_n)$ over $T$ to be a nodal divisorially $n$-marked
curve $C \to T$ and a {\em proper} morphism $f:C \to U$ such that, for each
geometric point $p$ of $T$, the group $\Aut_U(f:C_p
\to U, \Sigma_{1,p},\ldots,\Sigma_{n,p})$ is finite. Ideally, one should prove
the existence of the contraction of non-stable trees $C \to C'\to U$ for divisorially
$n$-marked maps $(C \to U,\Sigma_1,\ldots,\Sigma_n)$, and then deduce our
proposition using a presentation $R \double U$ of $\cM$. There are some
technical difficulties in carrying out such a local construction. We will use
the uniqueness of the local construction to deduce the uniquness in the
twisted case, but we will use a global approach for proving existence.
\end{remark}

\subsection{Twisted non-stable trees}
Before proving the proposition, we state two lemmas which illuminate the
geometric picture underlying the proposition. Consider a twisted curve $\cE$
over an algebraically closed field,
and a representable morphism $\cE \to \cM$ which is constant on moduli
spaces. If $p$ is the image point in $\bM$, then the reduction of its inverse
image in $\cM$ is the classifying stack $\cB G$ of the automorphism group of a
geometric point of $\cM$ lying over $p$. Thus the map $\cE \to \cM$ factors
through $\cB G$. 

\begin{lemma}\label{Lem:twisted-tree-1} 
Let $E$ be a tree of rational curves over an algebraically 
closed  
field $k$. Let $p\in E$ be a smooth closed point. Let $(\cE,\Sigma)$ be a
1-pointed twisted curve with associated coarse curve $(E,p)$. Let $G$ be a
finite group with order prime to $\chara k$, and let $\cE \to \cB G$ be a
representable morphism. Then $\cE \to E$ is an isomorphism, and the map $\cE \to
\cB G$ corresponds to the trivial principal bundle.
\end{lemma}
\proof We proceed by induction on the number of components of $E$, starting
from the case where $E$ is empty, when the statement is vacuous.  Let $D\subset
 E$ be a tail component; in case $E$ is reducible we may assume $p\notin D$,
 so in any case $D$ contains at most one special point $q$. Denote $D_0 := D
\setmin \{q\}$.  

Let $P \to \cE$ be the principal $G$-bundle corresponding to $\cE \to \cB
G$. Now restricting $P \to E$ to $D_0$, we obtain a tame  principal bundle
$P_{D_0} \to D_0$, which is  trivial, since
$\pi_1^{\rm tame}(D_0) = \{1\}$. 

Define $\cD$ to be the inverse image of $D$. The principal bundle $P_{\cD} \to
\cD$ has a section outside $q$, which automatically extends to $\cD$. Therefore
$P_{\cD} \to \cD$ is trivial, which means that $\cD \to \cB G$ factors through
$\Spec k \to \cB G$. Since this map representable, we have that $\cD$ is
representable.

By induction, the lemma holds over $\overline{E \setmin D}$. It follows that
$\cE$ is representable, and since $\pi_1(E)=\{1\}$, the principal bundle $P \to
\cE$ is trivial.
\qed

\begin{lemma}\label{Lem:twisted-tree-2}
 Let $E$ be a tree of rational curves over an algebraically closed 
field $k$. Let $p_1,p_2\in E$ be two distinct smooth closed points. Let
$(\cE,\Sigma_1,\Sigma_2 )$ be a
2-pointed twisted curve with associated coarse curve $(E,p_1,p_2)$. Let $G$ be a
finite group with order prime to $\chara k$, and let $\cE \to \cB G$ be a
representable morphism. Denote by $F\subset E$ the unique chain of components
connecting $p_1$ with $p_2$. Then 
\begin{enumerate}
\item $\cE \to E$ is an isomorphism away from the special points of $F$, namely
$p_1, p_2$ and the nodes  of $F$. 
\item There exists a cyclic subgroup $\Gamma\subset G$ such that the
image of the automorphism group of a geometric point $\Spec k\to E$ lying over
one of the special points of $F$ maps isomorphically to $\Gamma$.
\item\label{It:tree-cover-is-tree}
 Let $P \to \cE$ be the principal $G$-bundle associated to $\cE \to \cB
G$. Then each connected component of $P$ is a tree of rational curves.
\item For a geometric point $q \to \cE$ lying over a special point of $F$, we
let $\Gamma$ act on the tangent space of $\cE$ at $q$ via the isomorphism
above. Assuming $\cE$ is balanced, then the character of the action of $\Gamma$
on the tangent space to $\cE$ at $p_1$ is opposite to the character of the
action at $p_2$. 
\end{enumerate}
\end{lemma}
\proof Consider a connected component $E'\subset \overline{E\setmin F}$. Then
$E'$ is a tree of rational curves which is attached to $F$ at a unique point
$p'$. Let $\cE'\subset \cE$ be the reduction of the inverse image of $E'$ and
$\Sigma'$ be the reduction of the inverse image of $p'$ in $\cE'$. Then we can
apply the previous lemma to $\cE \to \cB G$ and conclude that $\cE' \to E'$ is
an isomorphism. Thus we may replace $E$ by $F$ and assume that $E=F$. 

Let $P \to \cE$ be the principal $G$-bundle associated with the morphism $\cE
\to \cB G$. Then $P \to E$ is a ramified covering which is tame and \'etale
outside of the special points. Let $Q \subset P$ be a connected component and
let $\Gamma$ be its stabilizer. Note that $\cE = [P/G]$, and since the action
of $G$ on the set of connected components of $P$ is transitive, we have that
$Q / \Gamma \to \cE$ is an isomorphism. Let $Q_1 \subset Q$ be an irreducible
component, and $E_1\subset E$ its image. Then, since the tame fundamental group
of $\PP^1\setmin \{0, \infty\}$ is cyclic, we know that $Q_1$ is a smooth
rational curve, and the map $Q_1 \to E_1$ is totally ramifies at exactly the
two special points of $E$ on $E_1$. This implies that the map of  dual
graphs  $\graph(Q)\to\graph(E)$ is \'etale, and since $\graph(E)$ is simply
connected 
this implies that this map of dual graphs is an isomorphism.  From this it
follows  that there is exactly one irreducible component of $Q$ over each
irreducible component of $E$ 
and the results are clear. \qed

\subsection{Proof of the proposition.} We may assume that $T$ is noetherian,
since 
$\cC\to T$ and $\cM\to T$ are of finite presentation.

{\sc Step 1: reduction to the representable case.} First we reduce to the case
where $\cC \to \cM$ 
is representable, as follows: let $R \double U$ be a presentation of $\cM$.
Consider the pullbacks $\cC_R \double \cC_U$ of $\cC \to \cM$ to $R$ and
$U$. These are divisorially $n$-marked twisted curves. Let $\cC_R \to C_R$ and
$\cC_U 
\to C_U$ be the respective coarse moduli 
spaces, which are divisorially $n$-marked  curves. There is an induced \'etale
groupoid $C_R \double C_U$ whose quotient 
is a twisted pointed curve $\overline\cC \to T$ with a representable morphism
$\overline\cC \to 
\cM$. The coarse moduli space of $\overline\cC$ is canonically isomorphic to
$C$. Note that, since $\cC\to T$ is tame, the formation of $C_R \double C_U$
commutes with base change on $T$, so the formation of $\overline\cC\to T$ also
commutes with base change on $T$.

Let us assume that $\cC\to T$ is balanced. In order to show that
$\overline\cC\to T$ is balanced as well, it is enough to consider the case
where $T$ is the spectrum of an algebraically closed field. Let $(V,\Gamma)$ be
a chart for  $\cC$. Without loss of generality 
we may assume $\Gamma$ is cyclic and has a unique fixed geometric point $q\in
V$, and the action of $\Gamma$ is free outside $q$. The morphism $\cC \to \cM$
induces a morphism $\Gamma \to \Aut(q\to 
\cM)$. Denote the kernel of this morphism by $H$ and the image by
$\overline\Gamma$. Then $(V/H,\overline\Gamma)$ is a chart for $ \overline\cC$,
which is evidently balanced. 

Thus from here on we will assume $\cC \to \cM$ is representable. 

{\sc Step 2: uniqueness.} Let $\cC \to \cC'\to \cM$ be a contraction as in the
proposition. Consider a presentation $R \double U$ of $\cM$ and consider fiber
diagram
$$ \begin{array}{ccccc}
  \cC_R & \double & \cC_U & \to & \cC \\
   \dar &         & \dar  &     & \dar\\
  \cC'_R& \double & \cC'_U& \to & \cC' \\
   \dar &         & \dar  &     & \dar\\
      R & \double &     U & \to & \cM. 
\end{array} $$

Note that $\cC_R$, $\cC_U$, $\cC'_R$,  and $\cC'_U$ are
representable. Moreover, it follows from Lemmas \ref{Lem:twisted-tree-1} and
\ref{Lem:twisted-tree-2} that the fibers of $\cC_U \to \cC_U'$ and $\cC_R \to
\cC_R'$ are precisely the non-stable trees in $\cC_U$ and $\cC_R$. 

   Lemma \ref{Lem:local-contraction-unique} below implies that
the groupoid $\cC'_R \double  \cC'_U$ is uniquely determined  as a local
contraction of non-stable trees in $\cC_R  \double  \cC_U$, and its formation
commutes with base change on $T$. Therefore $\cC \to \cC'\to \cM$ is unique,
and its formation also commutes with base change on $T$.

\begin{lemma}\label{Lem:local-contraction-unique}
Let $U \to T$ be a scheme of finite type,  $C\to T$ a divisorially marked
curve, and $C\to U$ a proper morphism. If $C\to C' \to U$ is a contraction of
the non-stable trees, then it is unique up to a canoncal isomorphism and its
formation commutes with base change.
\end{lemma}
{\bf Proof of lemma.} Note that the last statement follows from unicity.
Now $C'$ is uniquely determined as a topological space, since the fibers are
uniquely determined and the topology is the quotient topology since the map
$\pi:C \to C'$ is proper.
Therefore it is enough to check that the morphism  $\cO_{C'}\to \pi_*\cO_C$ is
an isomorphism.

Note that this morphism is injective since the kernel is supported in the locus
of special points of $C'\to T$, which contains none of the associated points
of $C'$.

 Denote by $Q$ the quotient sheaf: 
$$ 0 \to \cO_{C'}\to \pi_*\cO_C \to Q \to 0.$$

We need to show $Q=0$. We may assume that $T$ is the spectum of a local ring $R$. 
\begin{enumerate} 
\item Consider the case where $R$ is a field. We may pass to the algebraic
closure, and then the statement follows from the fact that nodal curves are
seminormal. 
\item Consider the case $\depth R> 0$. In this case the ideal of the support of
$Q$ has depth $>1$. This implies that the sequence $ 0 \to \cO_{C'}\to
\pi_*\cO_C \to Q \to 0$ splits locally on $C'$, by the $\cursext$
characterisation of depth. Let $\pi_*\cO_C \to \cO_{C'}$ be a local
splitting. It is enough to prove that this is injective. Let $f$ be a section
of the kernel of this splitting. By restriction to the central fiber we see
that $f \in \m\cO_C$, where $\m$ is the maximal ideal of $R$; moreover, as a
section of $\m\cO_C$, $f$ vanishes outside the non-stable trees. Since
$\m/\m^2\cO_C$ is a constant vector bundle on $C_0$, we have that the restriction
of $f \in \m\cO_C$ to the  central fiber is again zero. Inductively we show
that $f\in \m^k\cO_C$ for all $k>0$, and thus $f=0$. 
\item In the general case, denote by $I\subset \m$ the largest ideal of finite
length (so $R/I$ is either a field or of positive depth).
Let $J\subset I$ a an ideal of length 1. By induction on the length of $I$, we
may assume the statement holds true over $\tilde R := R/J$. 

The exact sequence 
$$ 0 \to J \to R \to  \tilde R \to 0$$ 
induces an exact sequence 
$$ 0 \to \cO_C \otimes J \to \cO_C \to \cO_{\tilde C} \to 0.$$
Note that $\cO_C \otimes J$ is isomorphic to $\cO_{C_0}$, where $C_0$ is the
central fiber. Applying $\pi_*$ we get an exact sequence
$$0 \to \pi_* \cO_{C_0} \to \pi_* \cO_C \to \pi_*\cO_{\tilde C} \to
R^1\pi_*\cO_{C_0}.$$

Note that the last term is zero. This follows by considering the inverse image
of an affine neighborhood of a  point on $C'$  containing the image of 
 a contracted tree $D\subset C_0$ and computing the cohomology
long exact sequence of 
$$0 \to \cO_{\overline{C_0\setmin D}}(-\Delta) \to \cO_{C_0} \to \cO_D \to 0.$$

Now we get a commutative diagram
$$\begin{array}{ccccccccc}
0 & \to & \cO_{C'_0}& \to &\cO_{C'}&\to&\cO_{\tilde C'} &\to& 0 \\
  &     &  \dar     &     & \dar   &   & \dar           &   & \\
0 & \to & \cO_{C_0} & \to &\cO_{C} &\to&\cO_{\tilde C}  &\to& 0
\end{array}$$
where the left and right columns are isomorphisms. This implies that the center
column is also an isomorphism, which concludes the proof
the lemma. \qed
\end{enumerate}

{\sc Step 3: construction of $\cC'$.} Consider the diagram
$$
\begin{diagram} 
\node{\cC}\arrow{e}\arrow[2]{s}{\pi}\arrow{se}\node{C}\arrow{s} \\
                                              \node[2]{C'}\arrow{s}\\
\node{\cM}\arrow{e}\node{\bM ,}\end{diagram}
$$
where $C \to C'\to \bM$ is the contraction constructed by Behrend--Manin.
Let $L'$ be a sufficiently relatively ample invertible sheaf on $C'/\bM$, and
let $\cL$ be its pullback to $\cC$. Here ``sufficiently ample'' means that the
degrees of $L'$ on the components of the fibers of $C'\to \bM$ are larger than
a certain integer $N$ to be specified later.

We define a quasicoherent sheaf of graded algebras 
$$\cA = \oplus_{r\geq 0}\pi_*\cL^{\otimes r}.$$

\begin{claim}
\begin{enumerate}
\item for any quasicoherent sheaf of $\cO_T$-modules $\cF$, we have $ R^1\pi_*
(\cL^{\otimes r}\otimes \cF)  = 0$ for all $r>0$.
\item For all $r>0$, the formation of $\pi_* \cL^{\otimes r}$ commutes with
base change on $T$, and the sheaf is flat over $T$.
\item The algebra $\cA$ is generated over $\cO_\cM$ in degree 1; in particular
it is locally  finitely generated.
\end{enumerate}
\end{claim}
\proof Note that any quasicoherent sheaf of $\cO_T$-modules $\cF$ we have
$R^2\pi_* 
(\cL^{\otimes r}\otimes \cF) = 0$, therefore the formation of $ R^1\pi_*
(\cL^{\otimes r}\otimes \cF)$
satisfies base change on $T$ (see \cite{G-EGA3}, 7.3.1).

 Let $U$ be an affine scheme and $U \to \cM$
a surjective \'etale morphism of finite type. Denote by $\cC_U$ the pullback of
$\cC$ to 
$U$. Let $D\subset \cC_U$ be the subcurve
consisting of components of fibers mapping to a point in $U$.
 Note that, since $T$ is noetherian, the number of topological types occuring
in the fibers of $D 
\to T$ is finite. 
Moreover, it follows from Lemma \ref{Lem:twisted-tree-1} and statement
\ref{It:tree-cover-is-tree} in Lemma 
\ref{Lem:twisted-tree-2} that the inverse images of the non-stable trees of
$\cC$ 
in $\cC_U$ are disjoint unions of trees of rational curves.

To check that $ R^1\pi_*
(\cL^{\otimes r}\otimes \cF)=0$, it is enough to check
the case where $\cF$ is the structure sheaf of a geometric point of $T$. This
is a consequence of the following elementary 
lemma, whose proof is left to the reader:

\begin{lemma}\label{Lem:pointwise-contracting-algebra}
Let $E$ be a quasi-projective nodal curve,  and let $D\subset E$ be the maximal proper
subcurve. Then there exists a number $N$, 
depending only on the topological type of $D$, such that the following holds:
Let  $\cL$ be an invertible sheaf satisfying  the following properties 
\begin{enumerate} 
\item for any component $D'\subset D$, the degree of $\cL$ on $D'$ is {\em
either} $>N$ or $0$, and
\item the union of the components $D'\subset D$ such that $\deg_{D'}\cL = 0$ is
a union of trees of rational curves. 
\end{enumerate}
Then 
\begin{itemize}
\item $H^1(E, \cL) = 0$, and 
\item the algebra $\oplus H^0(E, \cL^{\otimes k})$ is generated in degree 1.
\item Let $E' = \Proj\oplus H^0(E, \cL^{\otimes k})$. Then $E'$ is the
seminormal curve obtained by contracting the components of $D$.
\end{itemize}
\end{lemma}

This concludes the proof of part (1) of the claim. Also, once we prove part
(2), then the lemma above implies part (3).

For part (2), note that since  $ R^1\pi_*
(\cL^{\otimes r}\otimes \cF)  = 0$, the formation of $\pi_* \cL^{\otimes r}$ commutes with
base change on $T$. This means that for any quasicoherent sheaf
of $\cO_T$-modules $\cF$  we have the equality $\pi_* (\cL^{\otimes r}\otimes
\cF) = \pi_* 
\cL^{\otimes r}\otimes \cF$. But the functor sending $\cF$ to  $\pi_*
(\cL^{\otimes r}\otimes \cF)$ is evidently left exact, which means that the
sheaf is flat. 
This concludes the proof of the claim. \qed

 Define $\cC' := \Proj_{\cM}\cA \to \cM$.

{\sc Step 4: verification that $\cC' \to \cM$ is a twisted stable pointed map.}

 Since  $\cC'\to \cM$ is projective, it is representable.
Since $\pi_*\cL^{\otimes k}$ is flat over $T$ for all $k>0$, we have that $\cC'
\to T$ is flat. 

 Let $\Sigma_i^{\cC}$ be one of the
markings. We claim that the composition $\Sigma_i^{\cC}\to \cC \to \cC'$ is an
embedding, and we define $\Sigma_i^{\cC'}$ to be its image in $\cC'$. This, as
well as the
fact that $\cC'\to \cM$ satisfies statement (1) of the proposition, can be
checked on geometric fibers.  We claim that statement (2) can also be checked
on the level  of geometric fibers. First, the local criterion of flatness
implies that a proper morphism $C_1 \to C_2$ of flat $T$-schemes which is an
isomorphism on the geometric fibers of $T$ is an isomorphism. Also, the
formation of coarse moduli spaces of {\em tame} stacks commutes with arbitrary
base change (see Lemma \ref{Lem:tame-cms-pullback}). Since the moduli space
of a family of twisted nodal curves is flat, the claim follows.

So we assume that $T$ is the spectrum of an algebraically  closed field.

The formation of the algebra $\cA$ commutes with base change along \'etale
covers  $U \to \cM$. It follows from Lemma
\ref{Lem:pointwise-contracting-algebra} that over $U$, the curve $\cC'_U =
\Proj_U\cA_U$ is the 
contraction of the non-stable trees in $\cC_U$. Using descent and the uniqueness
of the contraction, we have that $\cC \to \cC' \to \cM$ satisfies statement (1)
in the proposition. Denote by $\overline{C'}$ the moduli space of $\cC'$. The
morphism $\overline{C'} \to C'$ is a bijection on geometric points and is an
isomorphism away from the special points, and since $C'$ is seminormal, the
morphism is an isomorphism. 
\qed

\appendix\section{Grothendieck's existence theorem for tame stacks}

\subsection{Statement of the existence theorem}
Consider a proper tame stack ${\cX}$ of finite type over a
noetherian complete  local ring $A$ with maximal ideal $\m$.
Then
${\cX}$ has a moduli space ${\bX}$, which is a proper
algebraic space over $A$.

For each
nonnegative integer  $n$ set ${\cX}_n = {\cX}\times_{\Spec A} \Spec
{A/ {\m^{n+1}}}$. We 
define the category of formal
coherent sheaves on ${\cX}$ in a rather simple-minded
way. A formal coherent sheaf
$\widehat{\cF}$ on ${\cX}$ is a collection ${\cF}_n$
of coherent sheaves on ${\cO}_{{\cX}_n}$, together with
isomorphisms ${\cF}_{n+1} \rest{{\cX}_n} \simeq {\cF}_n$ of sheaves
over ${\cO}_{X_n}$. A homomorphism 
$\widehat{\cF} \to \widehat{\cG}$ of formal coherent
sheaves as a compatible sequence of morphisms ${\cF}_n \to
{\cG}_n$, in the obvious way. We shall denote the category
of formal coherent sheaves on ${\cX}$ by $\fcoh({\cX})$.
From this description it is not even clear that $\fcoh({\cX})$ is
abelian. 

There is an obvious functor from the category $\coh( {\cX})$ of
coherent sheaves on ${\cX}$ to $\fcoh( {\cX})$, 
sending a coherent sheaf ${\cF}$ to the compatible system
of restrictions ${\cF}\rest{{\cX}_n}$.

\begin{theorem}\label{Th:groth-ex} If ${\cX}$ is proper and tame over
$\Spec A$, the functor from $\coh( {\cX})$ to $\fcoh(
{\cX})$ is an equivalence of categories.
\end{theorem}

This is of course well known if ${\cX}$ is an algebraic
space (see \cite{Knutson}, V, 6.3). The proof for stacks is not too
different from the one in \cite{Knutson}.

\subsection{Restatement in terms of $\widehat{\cO}_{\cX}$ modules} For the
proof, we will need a more manageable 
description of the category of formal coherent sheaves.
Consider the etale site $\xet$ of ${\cX}$, where the
objects are the \'etale morphisms $U \to {\cX}$, where $U$
is a scheme, and the arrows are the maps over ${\cX}$, and
coverings are defined as usual. On this site we have two
sheaves of rings; the usual structure sheaf, which we denote
by ${\cO}_{\cX}$, which sends each \'etale map $U \to
{\cX}$ to ${\cO}(U)$, and the completed structure
sheaf $\widehat{\cO}_{\cX}$, which sends an \'etale map
$U \to {\cX}$ to the completion of ${\cO}(U)$ with
respect to the inverse image of the maximal ideal of $A$. We
denote by $\widehat {\cX}$ the ringed site $\xet$, equipped
with the latter sheaf of rings $\widehat{\cO}_{\cX}$.

Notice that in case that ${\cX}$ is a scheme, the site
$\widehat{\cX}$ is not the \'etale site of the
corresponding formal scheme; however, the two categories of
sheaves of modules on $\widehat{\cX}$ and on the
\'etale site of the formal scheme are canonically equivalent,
so this does not really make a difference.

Consider a sheaf of $\widehat{\cO}_{\cX}$-modules ${\cF}$. Let $U \to
{\cX}$ be an \'etale 
map, where $U$ is a scheme, and let $\widehat{\cO}_U$ be
the restriction of
$\widehat{\cO}_{\cX}$ to the Zariski site of $U$. The
sheaf $\widehat{\cO}_{\cX}$ is actually supported on
$U_0$, so the restriction of ${\cF}\rest U$ to $U$ defines
a sheaf on the formal completion $\widehat U$ of the scheme
$U$ along
$U_0$. We say that ${\cF}$ is a coherent sheaf on
$\widehat{\cX}$ when ${\cF}\rest U$ is a coherent sheaf
on the formal scheme $\widehat U$ for all \'etale maps $U \to
{\cX}$. The category of coherent sheaves of modules over
$\widehat{\cX}$ will be denoted by $\coh(\widehat{\cX})$; it is
clearly an abelian category. 

There is an obvious functor $\coh(\widehat{\cX}) \to \fcoh( {\cX})$
that sends a coherent sheaf of 
modules ${\cF}$ over $\widehat{\cX}$ to the object
$\widehat{\cF} \eqdef ({\cF}\otimes_{\widehat{\cO}_{\cX}}
{\cO}_{{\cX}_n})$ of $\fcoh( {\cX})$. 
This functor is easily checked to be an equivalence of
categories; the inverse is obtained by taking inverse limits
in the usual fashion, i.e., one can send an object $({\cF}_i)$ of
$\fcoh( {\cX})$ into $\projlim {\cF}_i$, and 
this is a coherent sheaf of modules over $\widehat{\cX}$.

Now, the functor $\coh ({\cX}) \to \fcoh ({\cX})$ above
corresponds to the functor $\coh ({\cX}) \to
\coh(\widehat{\cX})$ which sends a coherent sheaf ${\cF}$ to $\widehat
{\cF} \eqdef {\cF}$. This induces a 
natural morphism $\opH^i({\cX}, {\cF}) \to
\opH^i(\widehat{\cX}, \widehat{\cF})$ from the cohomology
of
${\cF}$ on
${\cX}$ to the cohomology of the completion of ${\cF}$
on $\widehat{\cX}$. The first imoportant step of the proof
the theorem consists in proving the following form of
Zariski's theorem on formal functions.

\subsection{The theorem on formal functions}
In what follows we will use repeteadly the well known fact
that the cohomology of a coherent sheaf, considered as
a sheaf on the small \'etale site of a scheme, is the same
as the cohomology of the same sheaf on the Zariski site.

\begin{theorem}If ${\cX}$ is
proper over $\Spec A$, the morphism
$\opH^i({\cX}, {\cF}) \to
\opH^i(\widehat{\cX}, \widehat{\cF})$ is an isomorphism.
\end{theorem}

\begin{lemma} Let $({\cF}_i)$ be a formal sheaf on
${\cX}$, and set $\widehat{\cF} = \projlim
{\cF}_i$. Then the natural map $\opH^n( {\cX},
\widehat{\cF}) \to
\projlim\opH^n({\cX}, {\cF}_i)$ is an
isomorphism.
\end{lemma}

\proof This is very much the proof of the corresponding fact
in \cite{Knutson}, V, 2.19, eccept that it is, hopefully, correct. Let
$U \to {\cX}$ an surjective
\'etale map, where $U$ is an affine scheme. Then all the fiber
products $U^n = U \times_{\cX} U\times _{\cX} \ldots
U\times_{\cX}U$ are affine, because ${\cX}$ is
separated. The usual spectral sequence
$$
\cechh^p\big(U/ {\cX}, {\cH}^q({\cF}_i)\big)
\Longrightarrow \opH^{p+q}({\cX}, {\cF}_i)
$$
relating \cech{} cohomology and usual cohomology degenerates,
because the cohomology of ${\cF}_i$ over each $U^n$
degenerates, so cohomology is equal to \cech{} cohomology
over the covering $U \to {\cX}$. The same is true for the
projective limit
$$
\widehat{\cF} \eqdef \projlim {\cF}_i.
$$
For this it is enough to show that the cohomology of
$\widehat{\cF}$ on each $U^n$ is zero in positive degree.
This is the content of the lemma that follows.

\begin{lemma} The cohomology of the limit of a strict
projective system of coherent sheaves on the small \'etale
site of an affine noetherian scheme $X$ is zero in positive
degree.
\end{lemma}

\proof Let $X$ be the affine scheme, $({\cF}_i)$
the projective system, $\widehat{\cF}$ the limit. According
to a theorem of Artin (see \cite{Milne}, Theorem 2.17), the
\cech{} cohomology of
$\widehat{\cF}$ on the small \'etale site $X\et$ equals
its cohomology; so it is enough to show that, given an
\'etale map $U \to X$ of finite type, with $U$ affine,  the
\cech{} cohomology
$\cechh^\mini(U/X,\widehat {\cF})$ is acyclic in positive
degree. Let us denote by $K_\mini({\cG})$ the augmented
\cech{} complex
$$
\Gamma(X,{\cG}) \longrightarrow \Gamma(U, {\cG}) \longrightarrow
\Gamma(U^2, {\cG}) 
\longrightarrow \cdots
$$ of a sheaf ${\cG}$.

Now, each of the maps
${\cF}_{i+1} \to {\cF}_i$ is surjective, so ${\cF}_{i+1}(U^n) \to
{\cF}_i(U^n)$ is surjective for each 
$n$, because $U^n$ is affine and ${\cF}_i$ is coherent.
Thus the projective system
$$
\cdots \to K_\mini({\cF}_{i+1}) \to K_\mini({\cF}_i)
\to \cdots
$$
is a projective system of acyclic complexes satisfying the
Mittag-Leffler condition in each degree, and therefore its
projective limit
$$
\projlim K_\mini({\cF}_i) = K_\mini(\widehat{\cF})
$$
is acyclic.\endproof

To conclude the proof, it is enough to show that the map
$$
\cechh^n\left(U/{\cX}, \widehat{\cF}\right) \to \projlim
\cechh^n\left(U/{\cX}, {\cF}_i\right)
$$
is an isomorphism. The \cech{} complexes
$\cechc^\mini\left(U/{\cX}, {\cF}_i\right)$ form a projective
system of complexes which satisfy a Mittag--Leffler
condition in each degree. Furthermore the cohomology groups
$\cechh^n\left(U/{\cX}, {\cF}_i\right)$ are artinian
modules over $A$, so the corresponding projective systems
also satisfy a a Mittag--Leffler condition. The statement
follows from \cite{G-EGA3}, 0, 13.2.3 \endproof

Consider a moduli space $\pi\colon  {\cX}\to{\bX}$. 

\begin{lemma} The natural map $( \pi_*{\cF})_r \to
\pi_*({\cF}_r)$ is an isomorphism.
\end{lemma}

\proof This statement is local in the \'etale topology, so we
may assume that ${\cX} = [Y/ \Gamma]$, $Y = \Spec R$,
${\bX} =
\Spec R^ \Gamma$. The sheaf ${\cF}$ is associated with the
$\Gamma$-equivariant $R$-module $M$, and $\pi_*{\cF}$
corresponds with the
$R^ \Gamma$-module $M^ \Gamma$. Then the statement corresponds
to the fact that the natural map $M^\gamma / \m^{r+1}M^\Gamma
\to (M/ \m^{r+1}M)^\Gamma$ is an isomorphism. The
surjectivity follows immediately from the fact that any
invariant in the right hand side lifts to an invariant in
$M$. For the injectivity, take an element $m$ in $M^\Gamma$
which goes to 0 in $(M/ \m^{r+1}M)^\Gamma$; then we can write
$m = \sum_i a_im_i$, where $a_i \in \m^{r+1}$ and $m_i \in
M$. By averaging, we may assume that the $m_i$ are invariant,
so that $m \in \m^{r+1}M^\Gamma$. \endproof

Now, from the theorem on formal functions for algebraic
spaces (\cite{Knutson}, V, 3.1), we conclude that
$$
\opH^n( {\cX}, {\cF}) = \opH^n( {\bX}, \pi_*{\cF}) =
\projlim\opH^n( {\bX}, \pi_*{\cF}_r) = \projlim\opH^n( {\cX}, {\cF}_r). 
$$

\subsection{Algebraization of extensions, kernels, and  cokernels}
\begin{lemma} Let ${\cF}$ and ${\cG}$ be two
coherent sheaves on ${\cX}$. The natural map
$$
\cursext^n_{{\cO}_{\cX}}({\cF}, {\cG})^{\widehat{\
}}
\to \cursext^n_{\widehat{\cO}_{\cX}}(\widehat{\cF},
\widehat{\cG\,})
$$
is an isomorphism for all $n$.
\end{lemma}

\proof This is a local statement in the \'etale topology, so
it follows from the case that ${\cX}$ is a scheme.
\endproof

\begin{lemma}\label{Lem:ext-alg} Let ${\cF}$ and ${\cG}$ be two
coherent sheaves on ${\cX}$. The natural map
$$
\Ext^n_{{\cO}_{\cX}}({\cF}, {\cG})
\to \Ext^n_{\widehat{\cO}_{\cX}}(\widehat{\cF},
\widehat{\cG}\,)
$$
is an isomorphism for all $n$.
\end{lemma}

\proof This follows from the theorem on formal functions and
the lemma above, by considering the local global spectral
sequences for $\cursext^n_{{\cO}_{\cX}}({\cF}, {\cG})$, and
$\cursext^n_{\widehat{\cO}_{\cX}}(\widehat{\cF},
\widehat{\cG\,})$. \endproof

This implies that the functor $\coh( {\cX}) \to \coh(
\widehat{\cX})$ is fully faithful. We have to prove that
all formal coherent sheaves on $\widehat{\cX}$ are
algebraizable, that is, it is isomorphic to the completion of
a coherent sheaf on ${\cX}$.

\begin{lemma}\label{Lem:ker-coker-alg} The kernel and the cokernel of a
morphism  
of algebraizable formal sheaves are algebraizable.
\end{lemma}

\proof Let ${\cF}$ and ${\cG}$ be coherent ${\cO}_{\cX}$-modules,
$\widehat\phi\colon  \widehat{\cF}\to 
\widehat{\cG}$ a morphism of $\widehat{\cO}_{\cX}$-modules. By the
previous result, $\widehat\phi$ comes 
from a morphism of ${\cO}_{\cX}$-modules $\phi\colon 
{\cF}\to {\cG}$. Since completion is an exact functor,
the kernel and the cokernel of $\widehat\phi$ are the
completion of the kernel and cokernel of $\phi$. \endproof

\begin{lemma} Any extension of algebraizable formal sheaves
is algebraizable.
\end{lemma}

\proof If ${\cF}$ and ${\cG}$ are coherent ${\cO}_{\cX}$-modules, the
natural map 
$$
\Ext^1_{{\cO}_{\cX}}( {\cF}, {\cG}) \to
\Ext^1_{\widehat{\cO}_{\cX}}(\widehat{\cF}, \widehat{\cG}\,)
$$
is an isomorphism, by Lemma \ref{Lem:ext-alg}. \endproof

\subsection{Proof of the existence theorem}
 By
noetherian induction, in order to prove the existence theorem we may assume
that every coherent sheaf 
on ${\cX}$ is algebraizable when it is zero on an open
nonempty substack of ${\cX}$. Let $\rho\colon  Y \to {\cX}$ be a
surjective finite generically \'etale map, where 
$Y$ is a scheme; this exists by \cite{L-MB}, 16.6. If $\widehat Y$
is the completion of $Y$, then the diagram of ringed
Grothendieck topologies
$$
\begin{array}{ccc}
\widehat Y      & \stackrel{\widehat\rho}{\to} & \widehat{\cX}\\
\down           &                      &\down\\
Y               & \stackrel{\rho}{\to} &{\cX}
\end{array}
$$
is cartesian. From this we deduce the following Lemma.

\begin{lemma} Let ${\cF}$ be a coherent sheaf on $Y$.
Then the natural map $\widehat{\rho_*{\cF}} \to
\widehat\rho_*\widehat{\cF}$ is an isomorphism.
\end{lemma}

\begin{lemma} Let $\overline {\cF}$ be a coherent sheaf
on $\widehat Y$. Then the sheaf $\widehat\rho_*\overline {\cF}$ is a
coherent sheaf of ${\cO}_{\widehat\cX}$-modules.
\end{lemma}

\proof By restricting to an \'etale map $U \to {\cX}$,
where $U$ is an affine scheme considered with the Zariski
topology, this becomes obvious. \endproof

From these two Lemmas we see that a formal sheaf on
$\widehat{\cX}$ of the form $\widehat\rho_*\overline {\cF}$ is
algebraizable.

Now, we form the morphism $\sigma\colon  Y \times_{\cX} Y
\to {\cX}$. Let
$\overline {\cF}$ be a formal coherent sheaf on ${\cX}$; there are two
natural maps $\psi_1$ and $\psi_2$ 
from $\widehat\rho_*\widehat\rho^*\overline {\cF}$ to
$\widehat\sigma_*\widehat\sigma^*\overline {\cF}$ induced by the two
projections $Y \times_{\cX} Y \to Y$; let $\overline {\cK}$ be the
kernel of the diffence $\psi_1 - \psi_2$. The 
adjunction map $\overline {\cF} \to \widehat\rho_*\widehat\rho^*\overline
{\cF}$ factors through $\overline {\cK}$, and, by flat
descent, the induced map $\alpha\colon  \overline {\cF} \to
\overline {\cK}$ is an isomorphism over the open dense
substack of
${\cX}$ where the map $\rho$ is flat. By Lemma \ref{Lem:ker-coker-alg},
$\overline {\cK}$ is algebraizable; also, the kernel and
the cokernel of the map $\alpha$ are algebraizable, by the
induction hypothesis. So the image of $\alpha$ is
algebraizable, and therefore by the same Lemma $\overline
{\cF}$ is algebraizable. \endproof

\end{document}